\documentclass[12pt,reqno]{amsart}

\usepackage[dvips]{graphicx} 

\usepackage{amssymb, latexsym, amsmath, amscd, array,
}

\usepackage[all]{xy}

\usepackage{color}
\usepackage{amsthm}
\usepackage{latexsym}
\usepackage{amscd}
\usepackage{array}
\usepackage{verbatim}

\theoremstyle{definition}

\numberwithin{equation}{section}
\newcommand\N {{\mathbb N}} 
\newcommand\Q {{\mathbb Q}} 

\newcommand\R {{\mathbb R}}

\begin{document}
  
\title[Infinitesimals as an issue of neo-Kantian philosophy]
{Infinitesimals as an issue of neo-Kantian philosophy of science}

\author{Thomas Mormann}

\address{T. Mormann, Department of Logic and Philosophy of Science,
University of the Basque Country UPV/EHU, 20080 Donostia San
Sebastian, Spain} \email{ylxmomot@sf.ehu.es}

\author{Mikhail G. Katz}

\address{M. Katz, Department of Mathematics, Bar Ilan University,
Ramat Gan 52900 Israel} \email{katzmik@macs.biu.ac.il}

\begin{abstract}
We seek to elucidate the philosophical context in which one of the
most important conceptual transformations of modern mathematics took
place, namely the so-called revolution in rigor in infinitesimal
calculus and mathematical analysis.  Some of the protagonists of the
said revolution were Cauchy, Cantor, Dedekind, and Weierstrass.  The
dominant current of philosophy in Germany at the time was
neo-Kantianism.  Among its various currents, the Marburg school
(Cohen, Natorp, Cassirer, and others) was the one most interested in
matters scientific and mathematical.  Our main thesis is that Marburg
neo-Kantian philosophy formulated a sophisticated position towards the
problems raised by the concepts of limits and infinitesimals.  The
Marburg school neither clung to the traditional approach of logically
and metaphysically dubious infinitesimals, nor whiggishly subscribed
to the new orthodoxy of the ``great triumvirate" of Cantor, Dedekind,
and Weierstrass that declared infinitesimals \emph{conceptus nongrati}
in mathematical discourse.  Rather, following Cohen's lead, the
Marburg philosophers sought to clarify Leibniz's principle of
continuity, and to exploit it in making sense of infinitesimals and
related concepts.

\end{abstract}

%\date{\today}


\maketitle
 
\tableofcontents

\section{Introduction}
 
The %
%
%\footnote{issue in footnotes \ref{f23} }
%
traditional historical narrative concerning infinitesimals runs as
follows. The idea of infinitesimals has been with us since
antiquity. Mathematicians have used one or another variety of
infinitesimals or indivisibles without really understanding what they
were doing.  Eventually, infinitesimals fell into disrepute for
logical and philosophical reasons, as enunciated by Berkeley and
others.

Despite Berkeley's devastating criticism, mathematicians
continued to use them until the 19th century with more or less good
intellectual conscience. Finally (according to the traditional
narrative) Cauchy, followed by Cantor, Dedekind, and Weierstrass,
succeeded in formulating a rigorous foundation for the calculus in
terms of the epsilon-delta approach. Thereupon infinitesimals were
``officially" expelled from the realm of legitimate mathematics once
and for all. Or so it seemed.

This traditional narrative is, however, seriously incomplete.  Some 80
years after mathematics had allegedly dismissed ``infinitely small
magnitudes" and related concepts as pseudo-concepts once and for all,
in 1960 the mathematician Abraham Robinson claimed to have saved
infinitesmals from the bin of pseudo-concepts. Thereby, an improved
version of the traditional narrative goes, Robinson restored the
reputation of infinitesimals as legitimate mathematical entities by
means of his non-standard analysis.  From the perspectives of
mathematics and history of mathematics, such a completed version of
the traditional narrative is certainly to be considered as an
improvement.  Nevertheless it still suffers from serious shortcomings.
For centuries, the infinitesimal and related concepts were discussed
not only by mathematicians and scientists but also by philosophers
like Leibniz, Newton, Malebranche, and Berkeley. While the thought of
these classic 17th and 18th century authors has been extensively
studied by historians of philosophy and mathematics, not much is known
about the 19th century philosophical context in which the ``great
revolution in rigor" and the alleged dismissal of infinitesimals took
place. This is hardly acceptable.  A philosophically satisfying
account has to take into consideration the historical fact that the
``great revolution in rigor" in mathematical analysis, led by the
German mathematicians Cantor, Dedekind, and Weierstrass (henceforth,
CDW), took place when philosophy in Germany was dominated by various
currents of neo-Kantian philosophy.  Some of the neo-Kantian
philosophers had a keen interest in science and mathematics.  Indeed,
the issue of the infinitesimal was vigorously debated in neo-Kantian
quarters, as we will discuss in detail in this text.

Our main thesis is that the Marburg neo-Kantians elaborated a
philosophically sophisticated approach towards the problems raised by
the concepts of limits and infinitesimals.  They neither clung to the
obsolete traditional approach of logically and metaphysically dubious
infinitesimals,%
\footnote{\label{f6}Actually, there are good reasons to contend that
the infinitesimals of the traditional approach were neither logically
nor metaphysically dubious, though they were attacked as such by
George Berkeley.  Sherry~\cite{She87} dissected Berkeley's criticism
(of Leibnizian infinitesimal calculus) into its \emph{logical} and
\emph{metaphysical} components, and a closer inspection thereof
reveals that the Leibnizian system for differential calculus was both
more firmly grounded than the Berkeleyan criticism thereof, and free
of logical contradictions
(Katz \& Sherry 2012 \cite{KS2}),
(Katz \& Sherry 2013 \cite{KS1}),
(Sherry \& Katz 2013 \cite{SK}).}
nor whiggishly subscribed to the new orthodoxy of the ``great
triumvirate" (Cantor, Dedekind, Weierstrass) that insisted on the
elimination of infinitesimals from any respectable mathematical
discourse in favor of a new approach based on the epsilontic doctrine.
Instead, the Marburg school developed a complex array of
sophisticated, albeit not always crystal-clear, positions that sought
to make sense of {\em both\/} infinitesimals {\em and\/} limit
concepts.  With the hindsight enabled by Robinson's non-standard
analysis, the Marburg stance seems wiser than that of Russell, Carnap,
and Quine who unconditionally accepted the orthodox epsilontic
doctrine, along with its simplistic philosophical ramifications
stemming from a strawman characterisation of infinitesimals as a
pseudo-concept.

The outline of this paper is as follows. In Section~\ref{two}, we
recall the basics of the Marburg school's neo-Kantian philosophy of
science and mathematics.  In particular, we dwell on some crucial
features of the Marburg account that distinguishes it from the Kantian
orthodoxy.  We also recall the basic ingredients of Cassirer's
philosophy of science as a theory of the formation of scientific
concepts, pointing out the salient differences between `Aristotelian'
substantial concepts of common sense and the functional or relational
concepts of modern science. In Section~\ref{three}, we consider two
neo-Kantian attempts to elucidate the progressive conceptual evolution
in science with the help of mathematical metaphors, to wit, Natorp's
\emph{equational metaphor} and Cassirer's \emph{Cauchy metaphor}.  In
Section~~\ref{four}, we discuss several attempts by the members of the
Marburg school to elucidate some of Cohen's notoriously obscure theses
on the ``essence" of the infinitesimal.  The issue of the closing
Section~\ref{five} is Cassirer's road from infinitesimal to functional
concepts in his mature philosophy of science and mathematics.

\section{Neo-Kantian Philosophy of Science and Mathematics}
\label{two}

To set the stage, let us review some historical and philosophical
background. Following the collapse of German idealism after Hegel's
death in 1831, German philosophy once again returned to Kant
(cf.~Coffa 1991 \cite{Cof}).  Such a reorientation did not, however,
result in a new Kantian orthodoxy.  Rather, the emerging neo-Kantian
philosophy, subscribing to the maxim ``With Kant Beyond Kant" (Otto
Liebmann 1865 \cite{Lie65}) adopted some of Kant's ideas and at the
same time came to criticize the master.  The most important currents
of neo-Kantianism were the so-called Marburg school founded by Hermann
Cohen, and the Southwest or Baden school founded by Wilhelm Windelband
and Heinrich Rickert.  With some oversimplification
%
%\footnote{For instance, Carnap's doctoral supervisor Bruno Bauch
%belonged to the Baden school, although he mainly worked in the area of
%the philosophy of science.}
%
one may say that the Southwest school was mainly interested in matters
of {\em Geisteswissenschaften\/}, while the members of the Marburg
school were mainly engaged in the task of a philosophical
understanding of mathematics and the sciences.  We will therefore
concentrate on the neo-Kantian doctrines of the Marburg school and its
contributions to a philosophical understanding of the problems posed
by infinitesimals, limits, and related concepts.

Baldly characterizing an account in epistemology or philosophy of
science as {\em neo-Kantian\/} may suggest that such an approach is
rather similar to Kant's.  This would be an error.  All neo-Kantians
agreed that Kant's philosophy was a promising starting point for
modern epistemology and philosophy of science, but not a doctrine that
had to be followed literally. Not surprisingly, they vigorously
disagreed concerning the best way to go ``beyond Kant''. In this
paper, however, we will not dwell on the issue of whether or not the
Marburg neo-Kantian interpretation of Kant did true justice to Kant
(cf.~Friedman 2000 \cite{Fri00} and M.~K\"uhn 2010 \cite{Ku}), as our
main topic is \emph{neo-Kantian} rather than \emph{Kantian} philosophy
of science and mathematics.

An authoritative survey of the essence of the Marburg neo-Kantianism
was published by Natorp on the occasion of Cohen's seventieth birthday
in 1912, in the prestigious journal \emph{Kant-Studien}.  The article
\emph{Kant und die Marburger Schule} \cite{Na12} can be considered as
a kind of official position paper of the Marburg school.  More
precisely, Natorp sought to present the Marburg current as the true
heir, although not an epigone, of Kant's original philosophy. This
endeavor had two parts.  On the one hand, he emphasized the salient
differences among the Marburg school's neo-Kantianism, Kantian
orthodoxy, and rival contemporary philosophical currents such as the
Baden school and neo-Hegelianism; on the other hand, he pointed out
that the Marburg school preserved the true essence of Kant's
doctrines.

\subsection{The Transcendental Method}

Natorp emphasized that for the Marburg school, the true core of Kant's
philosophy was the \emph{transcendental method}
(\cite[p.~194f]{Na12}).%
\footnote{Whether or not the `transcendental method' \`a la Natorp was
also the core of Kant's philosophy need not concern us here.  Indeed,
many renowned Kant scholars deny this.}
Everything in Kant's system that did not fit well with this method had
to be given up by true Kantians. The transcendental method%
\footnote{\label{f3}The method was called ``transcendental'' since it
went beyond the cognition that is immediately directed onto the
objects.  A more detailed discussion appears in the main text at
footnote~\ref{f9}.}
deals with the problem of the possibility of scientific
experience. More precisely, pursuing the `transcendental method' as
the universal method of philosophy is contingent upon two
requirements:

\begin{quotation}
The first is a solid contact with the established facts of science,
ethics, arts, and religion.  Philosophy cannot breathe in empty space
of pure thought, where reason aims to fly high only on the wings of
speculative ideas. \ldots The place of philosophy \ldots is the
fertile lowlands of experience in a broad sense, i.e., it seeks to
take roots in the entire creative work of culture (science, politics,
art, religion) \ldots

The second, decisive requirement of the transcendental method is to
provide for these cultural facts (science, ethics, art, \ldots) their
conditions of possibility. In other words, philosophy, by following
the transcendental method, has to exhibit and to elaborate the lawful%
\footnote{The original German ``gesetzm\"assig'' is a key term of
neo-Kantian philosophy of science and difficult to translate.  Its
meanings range from \emph{order-generating} to \emph{exhibiting
regularity}.}
ground (``Gesetzesgrund'') or \emph{logos} of those creative acts of
culture (Natorp \cite[p.~197]{Na12}).
\end{quotation}

Restricting our attention to the cultural fact of science, we may
restate Natorp's thesis in more modern terms by saying that philosophy
of science in the transcendental mode has the task of rationally
reconstructing the evolution of scientific knowledge and the
conditions of its possibility.

The task of revealing the conditions of the possibility of scientific
experience binds the Marburg `transcendental method' to Kant's
original `transcendental logic', which, by definition, investigates
how it is possible that our concepts are related to real objects.
More precisely, transcendental logic in Kant's sense is concerned with
the origin, the content, and the limits of experiential knowledge.

Pursuing the transcendental method, the ``critical idealism'' of the
Marburg school is led to a \emph{genetic} epistemology and theory of
science that regards the ongoing process of scientific creativity as
its essential feature, more so than its temporary results.  Natorp put
it as follows: Knowledge is always in the state of ``becoming'', it is
never ``closed'' or ``finished''.  Something non-conceptually
``given'', in particular something allegedly intuitively ``given''
cannot be accepted as such.  A ``given'' is just another name for a
problem to be solved.%
\footnote{\label{gegeben}As a then popular neo-Kantian pun put it: An
object is not \emph{gegeben} but \emph{aufgegeben}.  This pun loses
its effectiveness in English: it claims that an object is not
``given'' (gegeben) but ``posed'' (aufgegeben) as a problem.  The
Russian equivalent \emph{dan/zadan} appears in a recent collection of
essays on neo-Kantianism, where it is mentioned in an essay by
Sebastian Luft as translated by N. Dmitrieva \cite[p.~121]{GD}, and in
an essay by T. B. Dlugach \cite[p.~224]{GD}.}
 
In other words, for the philosophers of the Marburg school the fact of
science was a ``fact of becoming'' (\emph{Werdefaktum}).  Accordingly,
the basic task of a truly ``Kantian'' philosophy of science was to
make explicit the method of science as ``the method of an infinite and
unending creative evolution of reason'' (Natorp \cite[p.~200]{Na12}).

The rejection of a non-conceptual given in any form brought the
Marburg neo-Kantianism in open conflict with one of the cornerstones
of Kant's epistemology, to wit, the dualism of conceptual
understanding and intuition, as a ``non-conceptually given''.  Indeed,
as the Marburg neo-Kantians argued against Kant's original position,
if one really follows the transcendental method in its true sense,
then
\begin{quotation}
it is virtually impossible, as Kant does, to maintain this dualism of
epistemic factors if one takes seriously the core idea of the
transcendental method. (\cite{Na12}, 201)
\end{quotation} 

For the Marburg neo-Kantians, in contrast to Kant, both the Kantian
categories and his forms of sensibility of space and time were purely
conceptual.  Kant's sharp separation of understanding and sensibility
as two complementary faculties of the mind had to be given up.

Since Kant's original transcendental logic as a logic of the
possibility of experience was closely related to the forms of
sensibility, the Marburg neo-Kantians were led to give up the
distinction between formal (deductive) logic and transcendental logic.
For them, there was only one logic, namely, the logic of the
transcendental method as the comprehensive logic of the conditions of
the possibility of scientific experience (cf. Heis 2010
\cite[p.~389]{Heis10}).

\subsection{Concepts in Mathematics and Science}
\label{rel}

The emphasis on the evolving character of scientific knowledge gave
the issue of the evolution of scientific concepts a central position
in the Marburg philosophy of science.  In particular, to Cassirer,
this meant that philosophy of science had to investigate the common
evolution of scientific \emph{and} mathematical concepts firmly
planted in the course of their historical development. From the
Marburg perspective, these two conceptual developments were two
aspects of the same problem.  In \emph{Substanzbegriff und
Funktionsbegriff} 1910 \cite{Cas10} (henceforth, \emph{SF}) Cassirer
wrote:

\begin{quotation}
[We should consider] physical concepts no longer for themselves but,
as it were, in their natural genealogy, in connection with the {\em
mathematical\/} concepts.  In fact, the physical concepts only carry
forward the process that is begun in the mathematical concepts, and
which here [in mathematics--the authors] gains full clarity.  The
meaning of the mathematical concept cannot be comprehended, as long as
we seek any sort of presentational correlate for it in the given; the
meaning only appears when we recognize the concept as the expression
of a \emph{pure relation}, upon which rests the unity and continuous
connection of the members of a manifold.  The function of a physical
concept also is first evident in this interpretation. (SF 1910
\cite[pp.~219-220]{Cas10}; p.~166 in the 1953 edition)
 \end{quotation}

Note that Cassirer employed terms such as ``continuous'',
``connection'', and ``manifold'' not only in their strict mathematical
sense, but also in a metaphorical sense.  Such usage of philosophical
and scientific concepts was typical of Cassirer's thought throughout
his entire philosophical career (see Orth 1996 \cite{Ort96}).
Following the Husserl scholar Eugen Fink, Orth refers to concepts used
in such a ``metaphorical'' way, as ``operative concepts'' as opposed
to ``thematic concepts''.  Operative concepts are concepts based on
``intellectual schemata'' without being fully explicit.  They serve as
``orientations in a conceptual field'' and are of a metaphorical
character.  As Orth points out, Cassirer had at his disposal a
particularly rich supply of operative concepts stemming from the
(Kantian) philosophical tradition and from contemporary science, in
particular physics and mathematics (cf. Orth, ibid. 111).
\emph{Continuity} and its relatives were among Cassirer's favored
operative concepts.

In contrast to many currents of contemporary philosophy of science,
the philosophers of the Marburg school regarded philosophy of the
sciences and philosophy of mathematics as being of the same ilk,
namely, as a study of \emph{relational} concepts.  It may be
considered as a pleasing confirmation of Cassirer's unified approach
to mathematical and physical concepts that one of the great
mathematicians of the 20th century, Hermann Weyl, subscribed to a
similar view, quite independently of Cassirer and apparently unaware
of the similarity of convictions.
 
Weyl sought to overcome the deficiencies of a purely formal conception
of mathematics such as Hilbert's without being forced to build
mathematics on a restricted base of a Brouwerian intuitionism.  Weyl
therefore proposed to seek help from physics.  More precisely, he
considered theoretical physics as the guiding example of a kind of
knowledge endowed with a meaning completely dfferent from that of the
common sense or phenomenal meaning.  Thus, in order to endow the
symbols of mathematics with a meaning, Weyl saw only one possibility:

\begin {quotation} 
\ldots [to] completely fuse mathematics with physics and assume that
the mathematical concepts of number, function, etc. (or Hilbert's
symbols) generally partake in the theoretical construction of reality
in the same way as the concepts of energy, gravitation, electron, etc.
(Weyl 1925 \cite[p.~30]{Wey25})
\end{quotation}

Weyl clung to this thesis of the essential similarity of
conceptualization in mathematics and the sciences till the end of his
life.  The conclusion he reached in his late essay \emph{A Half
Century of Mathematics} is fully in line with the neo-Kantian
approach:
\begin{quotation}
It is pretty clear that our theory of the physical world is not a
description of the phenomena as we perceive them, but [rather] is a
bold symbolic construction.  However, one may be surprised to learn
that even mathematics shares this character. (Weyl 1951
\cite[p.~553]{Wey51})
\end{quotation}

This agreement between Weyl and the Marburg school is all the more
remarkable since Weyl arrived at it from a rather different
philosophical background: he was influenced mainly by Husserl's
phenomenology, and never had shown in his mature age any affinity to
Kantian or neo-Kantian philosophy.

\subsection{Substantive versus Relational Concepts}
\label{23}

After these general remarks on the Marburg account of mathematical and
scientific concepts, it may be expedient to consider some concrete
examples in detail.  Thereby we can hope to illustrate what the
Marburg philosophers, in particular Cassirer, intended to convey by
their thesis of the relational (or functional) character of scientific
concepts.  Let us start with an elementary example.  At first view,
which is sometimes called the Aristotelian one, there appears to be a
close analogy between common-sense concepts such as \emph{rock} and
mathematical concepts such as \emph{number} or \emph{triangle}, in
that the concept \emph{rock} corresponds to the class of all empirical
entities that are rocks, i.e. the class of all entities that have all
the properties a rock is assumed to have; and similarly, the
mathematical concept \emph{number} is said to correspond to the class
of all mathematical ``entities'' that are numbers, i.e. the class of
entities that have all the properties numbers are assumed to have.

Cassirer rejected such an analogy.  According to him, the unity of
mathematical and scientific concepts was not to be found in any fixed
group of properties, but rather in the rules,%
\footnote{Mathematically speaking, the dichotomy is between
properties, i.e., unary relations, and binary (and higher) relations.}
which represented, in a lawful way, the mere diversity of objects that
``fall under the concepts" as their cases (i.e.,
instantiations). Elementary examples of relational scientific concepts
in this sense are mathematical formulas that describe arithmetic
series (i.e., sequences) such as~$1, 3, 6, 10, \ldots$.  For such a
series, the ``construction of unity'' is provided by a formula that
describes their generation according to some general law. For
instance, the series~$1,3, 6, 10, \ldots$ is characterized by the law
that the difference of the differences of its consecutive members is
always~$1$.  This fact is succinctly expressed by the
formula~$a(n)=n(n+1)/2,\;n\in\N$. The members of such a series do not
have a common property (in any ordinary sense of property) but appear
as cases of a common functional law.

More generally, Cassirer considered the formulas of mathematics,
physics, and chemistry as paradigmatic examples of relational
scientific concepts since they brought singular facts into a lawful
context.  Algebraic equations of geometric curves provide somewhat
less elementary examples.  Such equations can be used to describe the
movement of material bodies.  More precisely, they are conceptual
devices for embedding the individual perceived positions of a body in
a continuous, even, smooth trajectory. Continuity, smoothness and
other concepts of the infinitesimal calculus, are, however, highly
theoretical `ideal' concepts.  The embedding of singular data into a
continuous or smooth trajectory is anchored in a complex web of
crucial idealizing assumptions.

For the Marburg neo-Kantians, the indispensable role of idealizations
such as \emph{continuity} and \emph{smoothness} for modern science
demonstrated that the \emph{real} could be understood only through the
\emph{ideal}.  To Cohen, this translated into a statement that the
\emph{infinitesimal} was a core concept of any truly modern logic of
science.  Cohen appears to have assumed that the notion of the
infinitesimal necessarily underlies the concepts of continuity and
smoothness (which is technically speaking not the case from the
viewpoint of modern $\epsilon,\delta$ definitions of continuity and
smoothness).  In the opening chapter \emph{Infinitesimal-Analysis} of
his \emph{Logik der reinen Erkenntniss} he explicitly contended:
\begin{quote}
If logic is to be a logic of science, i.e., a logic of the
mathematized natural sciences, then it must be primarily the logic of
the principle of the infinitesimal.  If this is not the case and this
core principle does not occupy centre stage, then logic itself still
hasn't gained its proper centre, it still belongs to the past.  The
new scientific thought is that which since Galileo, Leibniz, and
Newton has become systematically efficient [and for which the
infinitesimal does play a fundamental role--the authors]. (Cohen 1902
\cite[p.~31]{Co02})
\end{quote}

According to Cohen, the Leibnizian \emph{principle of continuity} was
the key that had opened the gate toward such a truly modern logic of
the infinitesimal (see Section~\ref{four}).  Regrettably, however,
later generations of philosophers and scientists had not faithfully
followed Leibniz' lead.  Therefore, a logic of the infinitesimal was
still in its infancy.  It was incumbent upon the Marburg school to
develop the foundations of a working logic of the infinitesimal.

A paradigmatic example of a relational concept in physics was for
Cassirer the concept of \emph{energy}.  The utility of the concept of
energy is not to describe any new class of objects, alongside the
already known physical objects such as light and heat, electricity and
magnetism.  Rather, it signifies only an objective lawful correlation,
in which all these ``objects'' stand.  The meaning of the concept of
energy resides in the equations that it establishes among different
kinds of events and processes.  Energy in the sense of modern science
is not an object in the traditional sense, but a unifying perspective
that sheds light on a manifold of experiences.

This is rendered most evident by the functional identity of potential
and kinetic energy through which states are identified with temporal
processes:
 
\begin{quotation}
The two [moments of kinetic and potential energy] are ``the same'' not
because they share any objective property, but because they occur as
members of the same causal equation, and thus can be substituted for
each other from the standpoint of pure magnitude (SF 1910
\cite[p.~264-265]{Cas10}, (1953, p.~199))
\end{quotation}

\emph{Energy} cannot be understood as the conceptual counterpart of
something empirical out there.  Rather, it is to be understood as an
order-generating principle.  In this respect it resembles the notion
of number by which we make the sensuous manifold unitary and uniform
in conception (cf. SF (1910, 252), (1953, 189)).  In contrast to the
concept of number, the concept of energy is a genuine concept of the
empirical sciences.  Hence, since ``number'' and ``energy'' both
served as order-generating principles in essentially the same manner,
this was considered as another argument in favor of the Marburg thesis
that mathematics \emph{and} mathematized empirical sciences followed
the same rules of one and the same transcendental logic.  The concept
of energy shows that in modern science the allegedly objective
``things'' of common sense and traditional metaphysics are replaced by
a web of mathematically formulated relations that yield objectivity to
scientific knowledge.  Thereby the notorious Kantian
``things-in-themselves'' can be dispensed with:

\begin{quotation} 
We need, not the objectivity of absolute things, but rather the
objective determinateness of the method of experience. [SF, (1910,
428), (1953, 322)]
\end{quotation} 

Characterizing scientific knowledge by idealizing functional relations
reveals that it does not aim at a description of how the world
``really is''.  The concepts of modern science are not the mental
images of certain pre-existing objects; rather they are tools that
offer new unifying perspectives as he elaborated in full detail in his
magnum opus \emph{Philosophy of Symbolic Forms} (Cassirer 1923--1929,
1953--1957 \cite{Cas23}) (henceforth PSF), in particular in the third
volume that takes up many issues of SF.
 
Ideal gases, ideal fluids, etc. are not limiting cases approximated by
the more or less homogenous gases or the more or less ideal fluids
found in nature.  Rather, idealizing concepts such as perfect gases or
perfect fluids have an epistemological role.  They provide conceptual
perspectives that allow the formulation of general relational laws and
thereby they help to make sense of reality as a manifold of
experiences.  
 
Cassirer described this theoretical unification of the scattered data
of sensations as an embedding of an \emph{incomplete} empirical
manifold of sensations in a \emph{completed} conceptual manifold.
Typically, such embeddings can be carried out in a variety of ways.
In contrast to Kant, for the neo-Kantian Cassirer there were no fixed
forms that determined how this process was to be carried out.  Rather,
the ever-growing variety of conceptual completions of our experiences
is revealed in the historical evolution of science itself.  For
Cassirer, the paradigmatic example of such a conceptual completion was
Dedekind's completion of the rational numbers~$\Q$ to the real
numbers~$\R$.  The essential point of this completion was not that
some ``ideal'' numbers were ``added'' to the already existing rational
numbers, but that the relational system~$\R$ of real numbers provided
us with a new conceptual perspective to ``see'' more clearly the
conceptual essence of the rational numbers~$\Q$ themselves (cf. PSF
III, 392).

Although the processes of concept formation in mathematics and physics
are similar, they are not identical.  After all, there \emph{is} a
difference between mathematics and physics, and philosophy of science
has the task of elucidating this difference.  Roughly, Cassirer
considered conceptualization in mathematics as a simplified version of
conceptualization in physics:

\begin{quotation}
In contrast to the mathematical concept, however, in empirical science
the characteristic difference emerges that the construction which
within mathematics arrives at a fixed end,%
\footnote{\label{ff}Here Cassirer appears to express the view that in
mathematics the construction arrives at a fixed end.  This reading
seems hardly compatible with his general neo-Kantian outlook according
to which the essence of science resides in its unending evolution.  A
more plausible reading of Cassirer here would be that he contended
that mathematical concepts are relatively fixed with respect to
empirical concepts, just as the elements of an infinite convergent
series of numbers are fixed although the series itself may not reach a
fixed limit point in finitely many steps.}
remains in principle \emph{incompletable} within experience. But no
matter how many ``strata'' of relations we may superimpose on each
other, and however close we may come to all particular circumstances
of the real process, nevertheless there is always the possibility that
some relevant factor in the total result has not been calculated and
will only be discovered with the further progress of experimental
analysis (Cassirer 1910 \cite[p.~337]{Cas10}, (1953, 254)).
\end{quotation}

Factual and theoretical components of scientific knowledge cannot be
neatly separated.  In a scientific theory ``real" and ``non-real"
components are inextricably interwoven.  Not a single concept is
confronted with reality but a whole system of concepts.

\section{Two Guiding Metaphors of neo-Kantian philosophy of science}
\label{three}

In line with the essence of the \emph{transcendental method} that
conceived of science ``as an infinite and unending creative evolution
of reason'', the philosophers of the Marburg school considered
empirical or mathematical concepts (or theories as systems of
concepts) as stages in an ongoing process of an unending conceptual
approximation. Accordingly, the task of philosophy of science was to
investigate the conditions of possibility for such an evolution.
Remarkably, for this endeavor the Marburg neo-Kantians heavily relied
on that science whose conceptual evolution they sought to elucidate,
to wit, mathematics.  In other words, both Natorp and Cassirer, each
in his own way, sought to tap the resources of mathematics to
elucidate the structure of the conceptual evolution of the
sciences. To this end, they introduced certain mathematical metaphors,
by exploiting the mathematical concepts of approximation and
convergence for which the concept of the infinite played an essential
role.

\subsection{Natorp's \emph{Knowledge Equation}}
\label{31}

Perhaps the best-known of these metaphors in the hightime of
neo-Kantian philosophy was Natorp's ``equational metaphor'' that
compared the evolution of science with the solution of a numerical
equation.  According to it, coming to know an object
(``Erkenntnisobjekt'') was analogous to the process of solving a
numerical equation.  To be specific, the reader may keep in mind a
specific equation such as~$x^3 - x^2 + x - 1 = 0$.  

In line with Natorp's didactic intentions this equation has been
chosen to convey several ideas concerning knowledge and its objects.
First, the fact that it has several different solutions indicates that
the process of research may not lead to unique results.  Furthermore,
the fact that two of its solutions are imaginary reflects the fact
that the research process may lead to an expansion of the original
fundamental concepts one started with.  Note that the admission of
complex numbers as solutions transcends the conceptual space in which
the equation was formulated, since its coefficients are all integers.
What is still missing in this metaphor is the ``infinite character''
of the knowledge equation.  Natorp was aware of this shortcoming and
tried to remedy it (see below).
 
According to Natorp, the object of knowledge may be considered as an
``unknown~$x$ of the knowledge equation'':

\begin{quotation}
If the object is to be the~$x$ of the equation of inquiry, then it
must be possible to determine the meaning of this~$x$ by the nature of
this equation (i.e., the inquiry itself) in relation to its known
factors (our fundamental concepts). From this it must follow whether
and in what sense the solution of this problem is possible for us.
This is the very idea of the transcendental or critical method (Natorp
1903 \cite[p.~10]{Na03}).
\end{quotation}
 
Natorp added the following further elucidations.  The transcendental
method did not aim to extend our knowledge \emph{beyond} the limits of
the scientific method.  Rather, it sought to clarify the \emph{limits}
of scientific knowledge.  It was called ``transcendental'' since it
went beyond the cognition that is immediately directed onto the
objects, but aimed to obtain information about the general direction
of the path to be taken.%
\footnote{\label{f9}Here we elaborate further on the term
``transcendental'' as discussed in footnote~\ref{f3}.}
It did not provide us with any specific knowledge about an object
beyond experience.  Hence, following the established Kantian
terminology it was \emph{transcendental}, but not \emph{transcendent}.
 
Both Natorp and all his fellow philosophers of the Marburg school
viewed the object of knowledge, not as an unproblematic starting point
of the ongoing process of scientific investigation, but rather as its
limit.
%
%\footnote{\label{f30}I may have already raised this issue before
%(YES), but would it be helpful to replace ``limit'' by ``ultimate
%horizon''? NO. The NK always use the expressions here that have
%mathematical connotations.  Ultimate horizon is too poetic. OK, let's
%leave this discussion in so next time I am reminded not to raise the
%issue :-)}
%
This object was a problem to be solved.  In its various versions, this
equational account of knowledge can be found in virtually all of
Natorp's epistemological writings.  One might object that Natorp's
equational model of scientific cognition is far too simple in the
sense that the empirical objects hardly ever show up as solutions of a
finite equation such as the one considered above.  It seems hardly
plausible that physical entities such as ``proton'' or ``quark'' fit
in the conceptual framework of one physical theory without remainder.

Natorp did not ignore this difficulty and complemented his equational
account so as to counter this objection.  Elaborating the equational
model, he pointed out that the object of knowledge was not simply a
problem (``aufgegeben'') but an \emph{infinite} problem that could be
solved in finite time only approximately by finite creatures like
ourselves.  He thus sought to elude the trap of an overstated Hegelian
rationalism:

\begin{quotation}
Although we conceive of the object of knowledge ($= x$), similarly%
\footnote{A common criticism of the Marburg school's epistemology was
that it was dangerously close to an overstated Hegelian rationalism.
Hence Natorp, although he had to admit some similarity with Hegel, was
at pains to distance the Marburg neo-Kantianism from any sort of
Hegelian rationalism.}
as Hegel does, only in relation to the functions of knowledge itself,
and consider it \ldots as the~$x$ of the equation of knowledge, \ldots
we have understood that this ``equation'' is of such a kind that it
leads to an infinite calculation. This means that the~$x$ is never
fully determined by the parameters~$a, b, c \ldots$ of the equation.
Moreover, the sequence of the parameters \ldots is to be thought of as
being not ``closed'' but rather extendable further and further.
(Natorp 1912 \cite[pp.~211-212]{Na12})
\end{quotation}

\subsection{Cassirer's \emph{Convergent Series}}

As a second example of the usage of mathematical metaphors for the
elucidation of the meaning of philosophical ideas, let us now have a
closer look at how Cassirer conceptualized the guiding idea of the
Marburg school, to the effect that the evolution of scientific
knowledge could be understood as a conceptual approximation process.
In contrast to Natorp, Cassirer insisted that this approximation was
not assumed to converge to an \emph{externally} given limit, but
rather as defined by a general \emph{internal} rule, as is done in the
arithmetic of rational numbers.  To give a mathematical model for such
a conception of the evolution of science, we can characterize a
rational sequence such as~$3.1$,~$3.14$,~$3.141$,~$3.1415,\ldots$ as
being convergent, obviously without relying on an assumed existence of
an element in~$\Q$ to which it would converge.

Cassirer used this elementary mathematical insight to illustrate his
thesis that one may meaningfully speak of the convergence of
scientific theories without assuming that there is a fixed reality
``out there'' to which the sequence of our theories is expected to
converge.  To be specific, consider a sequence of positive numbers
such as~$(1/n)$ that converges to~$0$.  This may apparently suggest
that a convergent series%
\footnote{Cassirer used the term equivalent to ``series'' for the
mathematical entity usually referred to as a ``sequence''.  We have
retained Cassirer's terminology.}
of concepts (or theories) converges to some ultimate \emph{external}
entity (or ``reality''), just as the arithmetical series~$(1/n)$
converges to the real number~$0$ external to it, with 0 itself not
being a member of the series~$(1/n)$.  Cassirer vigorously rejected
such a realist ``exterior'' interpretation:
\begin{quotation}
The system (Zusammenhang) and the convergence of the series take the
place of an external standard of reality.  Both system and convergence
can be established and determined, analogously to arithmetic, entirely
by comparison of the serial members and by the general rule, which
they follow in their progress. (Cassirer 1910 \cite[p.~426]{Cas10},
(1953, 321))
\end{quotation}
As already 19th century mathematics had taught us, in order to be able
to speak meaningfully about a convergent sequence of numbers, it is
not necessary to assume that there ``really'' is a number to which the
sequence converges. Rather, an arithmetical series can be defined as
convergent if it satisfies an internal requirement that can be
formulated \emph{without} reference to a possibly inexistent external
limit.  Such an internal requirement is provided by Cauchy's
criterion.

Mathematicians have pushed this ``internalization'' of the concept of
convergence even further.  As is well known, the ``gappiness'' of the
rational numbers~$\Q$ (residing in the absence of limit points of
certain Cauchy-convergent sequences) may be overcome by {\em
completing\/}~$\Q$ in an appropiate way.  More precisely, one can
embed the rationals~$\Q$ into a set~$C(\Q)$ of appropriately defined
equivalence classes of Cauchy sequences.  Thereby it can be ensured
that in the new completed realm~$C(\Q)$ that englobes the rational
numbers~$\Q$ as a part, every Cauchy sequence has a limit point.  From
this viewpoint, a real number is an encapsulation, or reification,
%
%\footnote{\label{f34}Wow. These are terms borrowed from education
%theory.  I am not sure what they are doing here.  I probably put them
%in in the first place.  NO. YOU ALREADY WORRIED ABOUT THEM SOME MONTHS
%AGO. THEY ARE "MY" TERMS. THEY HAVE NOTHING TO DO WITH EDUCATION, BUT
%ARE METAPHYSICAL. THE INTENDED MEANING OF REIFICATION IS THAT A RATHER
%UNWIELDY PROCESS IS TREATED AS AN OBJECT THAT CAN BE HANDLED AS ONE
%OBJECT. For the purposes of this text it might be better to use
%``idealisation''.  Let me know what you think.  OK, let's leave this
%discussion in so that I am reminded nexxt time I want to raise the
%issue.}
%
of the concept of convergence of a Cauchy sequence.

Cassirer took these mathematical constructions to be more than mere
technicalities.  He considered them as the pattern for his internally
defined account of the continuous evolution of scientific knowledge:
\begin{quotation}
No single astronomical system, the Copernican any more than the
Ptolemaic, can be taken as the expression of the \emph{true} cosmic
order, but only the totality of these systems as they unfold
continuously according to a definite connection. (Cassirer 1910
\cite[p.~427]{Cas10}, (SF, 322))
\end{quotation}

In other words, for Cassirer the ``true cosmic order'' was not given
by a single theory but by a convergent series of theories. He did not
assert that our theories ontologically converge to a mind-independent
realm of substantial things as the substrate of a ``final'' theory.
His notion of theoretical convergence was epistemological rather than
ontological. He understood the approximation of theories essentially
as an epistemological progress within which the historical progression
of our theories continually approximates, but never reaches, any
ideally complete mathematical representation of the phenomena.  Such
an ideal representation is not waiting ``out there'' to be
approximated; rather, it resides in the reification of the
approximation process that comes into being through this very process
itself.  In other words, in the course of its history science
completes itself, so to speak, analogously as a convergent sequence
without limit gives rise to a corresponding sequence with limit point
in a suitably completed domain. For Cassirer, the paradigmatic example
of such a completion was Dedekind's completion of the rationals to the
real numbers.  But for him, the significance of Dedekind's
construction went beyond developing a more abstract version of
unending decimal expansions.  According to him, idealizing completions
were the essence of the modern empirical and mathematical sciences
(cf. Mormann 2008 \cite{Mor08}).

\subsection{Idealisations, completions, and infinitesimals}

The ``completion-friendly'' perspective of the Marburg school on the
conceptual evolution of science had important consequences for matters
infinitesimal.  With respect to completions of number systems by
infinitesimals, the thesis of the ``incompletability" of the
conceptual evolution in science and mathematics suggested that
Cassirer's account (and that of Marburg neo-Kantianism in general) had
no built-in source of resistance to further ontological extensions
beyond the rational and the real numbers.  On the contrary, according
to its own rules, the neo-Kantian approach would have welcomed the
advent of the hyperreals \`a la Edwin Hewitt \cite{Hew48} and Abraham
Robinson \cite{Ro61, Ro66} and related developments.  As we shall see
in the following sections, the insistence on the openness of the
evolution of mathematical and physical concepts brought the Marburg
philosophers in conflict with mathematicians and logicians such as
Cantor and Russell, who adhered to more realist accounts that
considered concepts as more or less direct descriptions of what there
\emph{is} (in the empirical or in an ideal realm), instead of
conceiving of them as epistemological tools for progressing in the
task of making sense of some aspects of the world.  According to
Cassirer, both empiricism \`a la Bacon and Berkeley, and naive
realism%
\footnote{\label{ff13}``Realism'' is taken here as an unsaturated term,
i.e., questions of realism arise with respect to a certain subject
matter, e.g., realism with respect to atoms, values, mathematical
objects, possible worlds, causality, or macroscopic material objects.
Cantor and Russell were realists with respect to mathematical
entities.  Other names for this sort of realism are ``platonism'',
``platonist realism'', or even ``platonist idealism''. Cassirer, a
self-proclaimed ``critical idealist'' was not a partisan of platonist
idealism.  On the contrary, he vigorously criticized it.}
with respect to ``ideal'' mathematical entities, \`a la Cantor and
Russell, must be rejected:

\begin{quotation}
For the existence of the ideal, which can alone be critically affirmed
and advocated, means nothing more than the objective logical necessity
of idealization.  (SF 1910 \cite[p.~170]{Cas10}, (1953, 129))
\end{quotation}
And again:
\begin{quotation} 
The relation beween the theoretical and factual elements at the basis
of physics cannot be described in this simple way. It is a much more
complex relation, it is a peculiar interweaving and mutual
interpenetration of these two elements, that prevails in the actual
structure of science and calls for clearer expression logically of the
relation between principle and fact. (Cassirer (1910, 172), (1953,
130))
\end{quotation}

According to Cassirer, no scientific theory directly relates to the
facts of perception.  Rather, such a theory relates to the ideal
limits, which we substitute for the facts of perception. Thus, we
investigate the impact of bodies by regarding the masses, which affect
each other, as perfectly elastic or inelastic.  We study perfect
fluids even though no such are to be found.  In other words, Cassirer
sought to present his ``Critical Idealism'' as a theoretical framework
that overcame both a naive empiricism and a misled platonist
idealism. 
 
For the Marburg neo-Kantians, who had always emphasized the essential
unity of mathematics and empirical science, the new relational logic
(Frege, Peano, Russell and others; see e.g., Gillies \cite{Gi92}) was
part and parcel of a single comprehensive transcendental logic of
science that was emerging in the course of the history of science.
For them, it was a fundamental philosophical mistake of the `logicism'
of Frege and Russell to ignore its transcendental character in
conceiving of it as a purely formal device neatly separated from the
empirical realm.  For the neo-Kantians, this was just another example
of a philosophically untenable dualism, analogous to the Kantian
dualism between the conceptual and the sensual.%
\footnote{\label{f40b}See Section \ref{two}.}
%
%commented out material left as a historical record, to eb ignored by
%both authors: IT MAY BE HELPFUL TO ELABORATE A BIT AS TO WHY EXACTLY
%THE READER IS INVITED TO VISIT SECTION \ref{two} AND WHAT RELEVANT
%MATERIAL IS DISCUSSED THERE.

One of Cassirer's criticisms of Russell's philosophy of logic was that
it succumbed to a naive ``platonizing idealism" (see also J.~Heis
\cite[p.~386]{Heis10}) since it insisted on a strict separation
between the logico-mathematical conceptual realm, on the one hand, and
the empirical realm, on the other.  According to Cassirer, this stance
expressed itself in an outdated dualistic metaphysics that was bound
to lead into unsolvable, self-inflicted pseudo-problems (cf. SF 1910
\cite[pp.~313, 359]{Cas10}, (1959, pp. 237, 271)).  For instance, the
applicability of mathematics in the mathematized natural sciences
became an unfathomable mystery.%
\footnote{Perhaps Wigner put forward the most influential plea for this attitude in terms of ``the unreasonable effectiveness of
mathematics in the natural sciences" (cf. \cite{Wi60} and
\cite{Gra08}).}

\section{Three attempts to  make sense of Cohen}
\label{four}

The philosophers of Marburg neo-Kantianism considered themselves
consciously as members of a well-defined school under the leadership
of Hermann Cohen as the school's founder.  Even when they held
differing opinions concerning a philosophical question, they sought to
minimize their differences vis-\`a-vis outsiders.  An important
example of such ``school discipline'' concerns the concept of
infinitesimals and how it was dealt with in the treatises of the
school's leader Cohen.

Even for sympathetic readers, it is often difficult to make sense of
Cohen's writings. Such a difficulty is not limited to contemporary
readers accustomed to doing philosophy in a more analytic style.
Already the young Husserl in 1886 complained in a letter to Brentano
that Cohen's allegedly ``scientific philosophy'' was actually nothing
but ``nonsensical profundity" or ``profound nonsense" (see Mohanty
2008 \cite[p.~3]{Moh}).  As evidence for this strong claim Husserl
took Cohen's theory of the principle of continuity; see \emph{Das
Prinzip der Infinitesimalmethode und seine Geschichte (The Principle
of the infinitesimal method and its history)} (henceforth
\emph{Prinzip}) 1883 \cite[\S 40ff]{Co83}.

We argue that Husserl's radical verdict was not entirely justified.
For this purpose we will rely on three different attempts to make
sense of Cohen's ``scientific philosophy'' undertaken by three members
of the Marburg school, namely Ernst Cassirer, Dimitry Gawronsky, and
Paul Natorp.  The works involved are Cassirer's \emph{Leibniz' System}
(1902), \emph{Substanzbegriff and Funktionsbegriff, SF} (1910),
Gawronsky's \emph{Das Urteil der Realit\"at} (1910), and Natorp's
\emph{Die logischen Grundlagen der exakten Wissenschaft} (1910).  As
we will see, it reveals that the evolving thought of the Marburg
school about infinitesimals and related concepts was not monolithic
and without cracks.  This is not surprising, since Cohen's
\emph{Prinzip} and Cassirer's SF are separated in time by almost
thirty years, not to mention the rather different technical styles and
the scientific background of their authors.  Thus, Cohen's education
in logic did not correspond to the state of the art at the beginning
of the 20th century.  He apparently never took proper notice of Frege,
Russell, or any other contemporary logician. On the other hand,
Cassirer, Gawronsky, and Natorp were aware at least partially of the
new developments in logic and mathematics, and sought to adapt the
Marburg school's philosophical stance to the new circumstances
(cf.~Cassirer 1907 \cite{Cas07}).

\subsection{The Point of departure: Cohen's \emph{Prinzip}}

Cohen treated the issue of the infinitesimal first in \emph{Das
Prinzip der Infinitesimalmethode und seine Geschichte} and later, in
\emph{Logik der reinen Erkenntnis}.  The \emph{Logik} may be
considered as a continuation and philosophical elaboration of
\emph{Prinzip}.  The most significant difference is that the last
traces of any sort of Kantian `intuition' in the constitution of
infinitesimals are eliminated.  \emph{Pure thought},%
\footnote{\label{f44}The expression \emph{pure thought}
(``\emph{reines Denken}") was the technical term designed by the
neo-Kantians as the successor concept of the two separate components
of Kantian epistemology, to wit, ``concepts" and ``intuitions".
%
%I.e., ``transcendental logic''? RATHER NOT. I'd prefer the expression
%"pure thought", since this was the standard NK term for the union of
%the two Kantian components of "concepts" and intuitions".
}
and pure thought alone, takes care of matters infinitesimal.  Cohen
even went so far as to contend that the infinitesimal was to be
considered as the most important and most typical issue of pure
thought, \emph{\"uberhaupt}.%
\footnote{\label{f45}According to Cohen transcendental logic was the
logic of the infinitesimal (see Subsection~\ref{23}).  Infinitesimals
provided the key example of how the real was clarified by the ideal.}
In his \emph{Logik}, Cohen appears to assume that the reader has read
and digested the argumentation of \emph{Prinzip}.  That is to say,
mathematically there is nothing new in \emph{Logik} that cannot be
found already (usually more fully elaborated) in \emph{Prinzip}.  The
main purpose of \emph{Logik} is rather to explicate Cohen's account of
the philosophical or metaphysical presuppositions and ramifications of
\emph{pure thought}, centering on the notion of the infinitesimal.

In \emph{Prinzip}, Cohen still sought to establish the existence of
infinitesimals with the aid of an intuition in the sense of Kant,
while in \emph{Logik} any vestige of Kantian intuition had vanished.
According to mature Marburg neo-Kantianism, intuitions played no role
in scientific knowledge.  This must not be misunderstood: Cohen's
target was \emph{Kantian} intuition, not intuition in the everyday
understanding of the term. Recall that Kant distinguished between two
basic components of human representations, to wit, \emph{concepts} and
\emph{intuitions}. These correspond to two essentially different
cognitive faculties: understanding and sensibility.  Only if these
faculties are united can knowledge be achieved.  The neo-Kantians
argued against this two-tiered Kantian epistemology, not against the
common-sense claim that for the individual scientist, ``intuitions"
may play an important role in the context of discovery.

Though this may not be crystal-clear already from the first pages of
\emph{Prinzip}, in later works, in particular in \emph{Logik}, Cohen
repeatedly emphasized this difference between orthodox Kantianism and
Marburg neo-Kantianism.  The rejection of the significance of
``Kantian intuition" for scientific knowledge was, as is well-known, a
common trait of all members of the school.  Hence, critizising Cohen
for relating the concept of the infinitesimal to some sort of Kantian
intuition would be a gross misunderstanding.  Another more subtle,
albeit quite common misconception would be to ascribe to Cohen the
claim that he sought to locate the ``problem of the infinitesimal" in
the realm of a psychologistically conceived epistemology
(\emph{Erkenntnistheorie}).

Already in \emph{Prinzip}, Cohen had insisted that the problem of the
infinitesimal could be properly treated only in what he referred to as
\emph{the logic of science} a.k.a. ~\emph{transcendental logic}.
According to him ``the logic of science must be the logic of the
principle of infinitesimal calculus'' (cf.~\emph{Logik}, already
quoted in Section~\ref{23}).  He was eager to point out that a glance
at the literature revealed that the logic of his day had not yet
recognized the decisive logical significance of the infinitesimal
principle. In other words, and to give it a more personal twist, in
\emph{Logik} he admitted that his \emph{Prinzip} had not yet found the
recognition it deserved.

This situation did not change in the ensuing decades.  Thus, in the
otherwise rather comprehensive survey by Haaparanta, \emph{The
Relation between Logic and Philosophy, 1874 -- 1931} (Haaparanta 2009
\cite{Haa}), Cohen's and, more generally, the Marburg account of logic
as transcendental logic of science was completely ignored.

Cohen blamed Kant, to some extent, for this myopic conception of
logic. According to him, Kant had failed properly to understand the
role of the infinitesimal for a true critique of pure reason. Instead,
he introduced pure sensitivity as a second ingredient of knowledge,
whereby the independence of pure thought had been compromised
(cf. Cohen (1902) p.~32)).

Let us now take a closer look at Cohen's general conception of logic.
Cohen's conception incorporates both modern and obsolete ideas in a
peculiar mixture.  The first thing to note is that for Cohen logic is
independent of anything else: logic is neither a branch of psychology
nor a branch of linguistics.  The laws of logic, i.e., the laws of
``pure thought" are neither psychological laws nor grammatical laws.
Rather, according to Cohen, the laws of logic constitute the core of
pure thought. Thoughts in Cohen's sense should not be confused with
``sensations" (\emph{Empfindungen}) or individual mental
representations (\emph{Vorstellungen}).  Cohen vigorously rejected any
dependence of logic on another realm.  For him, logic described the
activity of pure thought that took place in the ongoing conceptual
evolution of science: ``The thought of logic is the thought of
science.  Thought constitutes the foundation of being" (Cohen 1902
\cite[pp.~17, 18]{Co02}).
 
Cohen's conception of logic is a far cry from any modern post-Fregean
or post-Russellian conception of logic.  There are no axioms,
inferential rules or anything of that sort.  Rather, his logic
follows, at least superficially, quite closely the patterns of a
Kantian (or even pre-Kantian) \emph{Urteilslogik} (judgment logic).
Thereby Cohen's \emph{Logik} reveals a somewhat paradoxical relation
to Kant.  On the one hand, it can be read as a definitive parting of
ways with Kantian orthodoxy, in particular by giving up the basic
structure of the Kantian philosophical system, namely the distinction
between the two pillars of pure logic and pure sensitivity.  On the
other hand, Cohen formally followed the architectonics of the
\emph{transcendental logic} of the \emph{Critique of Pure Reason} when
he mimicked in his \emph{Logik} Kant's system of categories and
judgments.  In close analogy to Kant's pair of $4 \times 3$ schemata
of categories and judgments, Cohen set up the following~$4 \times 3$
schema distinguishing between four classes of judgments each
consisting of three types of judgments:
\begin{itemize}
\item
The judgments of laws of thought\\ (Origin, Identity, Contradiction)
\item
The judgments of mathematics\\ (Reality, Majority, Totality)
\item
The judgments of mathematized natural sciences\\(Substance, Law,
Concept)
\item
The judgments of methodology\\(Possibility, Contingency, Necessity)
\end{itemize}
 
For the purposes of the present text, it is unnecessary to dwell upon
Cohen's schema in full detail.  But the following remarks may be in
order.  Superficially, Cohen's and Kant's tables of judgments are
quite similar.  Both exhibit the~$4 \times 3$ schema, and, moreover,
Cohen's ``judgments of methodology" and virtually identical with
Kant's ``judgments of modality".  Behind this formal similarity,
however, are lurking deep conceptual differences.  The first is that
Kant's schema is based on what he called general logic, based on the
standard formal logic of his time.  In contrast, Cohen's table is
deeply soaked with contentious assumptions of his account of
transcendental logic. Moreover, it is not concerned with knowledge in
general, but with knowledge of the mathematized natural sciences.  In
this respect, Cohen's epistemological perspective was considerably
narrower than Kant's.  Perhaps even more important is the difference
between Kant and Cohen concerning the problem of how to conceive of
the relation between the table of judgments and the corresponding
table of categories.  In Cohen's \em{Logik}\em, there is nothing that
even remotely resembles Kant's famous transcendental deduction of the
categories.

Kant invested immense efforts in the task of this deduction that
resulted in a strict and rigid 1-1 correspondence between the~$4
\times 3$ items of the table of judgments and the~$4 \times 3$ items
of the resulting a priori categories.  In contrast, Cohen was content
with the vague assertion that there was a mutual ``correlation''
between judgments and categories relying on the bland metaphorical
explication that ``the category is the aim of the judgment, and the
judgment is the road to the category" (\emph{Logik}, p.~47).  In
particular he gave up the 1-1 correspondence between categories and
judgments, and allowed that every category might be contained in
several judgments and every judgment might be contained in several
categories.  Actually, this vagueness and indecision is no
co-incidence.  Kant's categories were intended to be valid \emph{a
priori}, once and for all.  There was no change or evolution in the
categorical schema that Kant had set up in the \emph{Kritik}.  In
contrast, Cohen repeatedly emphasized the evolving character
(\emph{``Werdecharakter"}) of scientific knowledge, for instance, when
he contended that the ``truly creative elements of scientific thought
reveal themselves in the history of scientific thought''
(\emph{Logik}, 46), but, obviously, this dynamic character of
scientific knowledge was hardly compatible with a Kantian schema of
fixed \emph{a priori} categories.

For contemporary philosophy of science, Cohen's half-baked proposal of
how to reconcile the categorical structure of scientific knowledge and
its historical character can be of historical interest at best.%
 
Nevertheless one may note that Cohen's problem, as we may call it, has
remained on the agenda of virtually all accounts of philosophy of
science that have been inspired by Kant up to this very day, as
exemplified, for instance, by Reichenbach's reformulation of the
Kantian \emph{a priori} in the 1920s, up to Michael Friedman's
neo-neo-Kantian \emph{Dynamics of Reason} (Friedman 1999 \cite{Fri99})
where the author seeks to reconcile historical and the \emph{a priori}
aspects of scientific knowledge as a synthesis of ideas taken from
Kuhn, Cassirer, and Carnap.

As was already mentioned at the beginning of this section, we do not
aim at an exhaustive treatment of Cohen's logic of judgments.
Instead, we aim to shed some light on Cohen's often obscure analyses
by consulting the writings of other members of the Marburg school, to
wit, Ernst Cassirer, Dimitry Gawronsky, and Paul Natorp.  These
authors may be helpful in elucidating their master's thoughts, as all
of them intended not to deviate from Cohen's ideas unless absolutely
necessary. This does not mean, of course, that they actually provided
faithful and accurate interpretations of Cohen's account.
Nonetheless, the works of these three philosophers can be read as
sympathetic readings that try to make the best out of Cohen.

\subsection{Cassirer's \emph{Leibniz' System}}
\label{ca}

Cassirer's first philosophical works were his 1899 dissertation
\emph{Descartes' Kritik der Mathematischen und Naturwissenschaftlichen
Erkenntnis}%
\footnote{\label{f20}\emph{Descartes' Kritik} was published as the
first part of \emph{Leibniz' System}.}
and his \emph{Leibniz' System in seinen wissenschaftlichen Grundlagen}
(Cassirer 1902 \cite{Cas02}).
 
Cassirer eventually published his Descartes and Leibniz texts together
as one book proposing that Descartes may be conceived of as a
forerunner of Leibniz.  More precisely, according to Cassirer,
Descartes and Leibniz may be considered as two stations of the long
and winding road toward an idealistic conception of science.
Provisionally, this was achieved in Kant's account; after Kant it
found its contemporary expression in the philosophy of science, or the
scientific philosophy, of the Marburg school.  Cassirer's
\emph{Leibniz' System} is a 400 page long text.  In its ten chapters
Cassirer seeks to treat Leibniz' philosophical and scientific
achievements within logic, mathematics, mechanics and metaphysics.
The latter is understood in a broad sense, including Leibniz'
reflections on issues such as ``the problem of consciousness", ``the
problem of the individual", and ``the concept of the individual in the
system of \emph{Geisteswissenschaften}.''  By far the largest (and for
our purposes most interesting) chapter is the fourth, dealing with The
Problem of Continuity.  It comprises not less than 70 pages.  In this
section we will mainly concentrate on this chapter of \emph{Leibniz'
System}.

\emph{Leibniz' System} is engaged in the ambitious task of presenting
the conceptual evolution of modern science by presenting the
achievements of the two geniuses of Descartes and Leibniz.  The
starting point was Descartes' overcoming of the medieval conception of
science.  According to Cassirer's fundamental thesis, the
\emph{philosophical} systems of these men could not be understood by
separating them from their \emph{scientific} achievements, to wit,
analytical geometry in the case of Descartes, and the infinitesimal
calculus in the case of Leibniz:
 
\begin{quotation}
Through the discovery of analytic geometry Descartes lays the
foundations for the modern way of scientific thinking, which finds its
mature expression in the infinitesimal calculus.  \ldots The synthesis
of philosophy and science which is carried out thereby must not be
conceived of as a mere juxtaposition. \ldots One must attempt to
identify a common basis of these thoughts (Cassirer 1902
\cite[pp.~1-2]{Cas02}).
\end{quotation}
Identifying such a common base will lead to a more profound
understanding of the role of Descartes' system in the historical
evolution of critical epistemological idealism and its continuation
and completion in Leibniz and Kant and, one may add, in the idealism
of the Marburg school.  So much for Descartes.  In the remainder of
this section, we will concentrate on Cassirer's study of Leibniz as
one of the most important early sources for the philosophy of science
of the Marburg school.

According to Cassirer, for Leibniz mathematics was primarily an
instrument of scientific research and a presupposition for the
discovery of a new concept of nature, rather than an aim in itself
(Cassirer 1902, 99).  This became fully evident through the ``new
mathematics", i.e., the infinitesimal calculus.  Following Cohen,
Cassirer contended that for Leibniz' thought the concept of the
infinitesimal was to be considered of fundamental importance, not only
with respect to mathematics, but much more generally, also for
Leibniz' philosophical understanding of the mathematized empirical
sciences, a new concept of nature in general, and his metaphysics in
general. Indeed, \emph{Leibniz' System} is to be considered only as
Cassirer's first attempt to contribute to this overall programme of a
genuine Marburg philosophy of science inaugurated by Cohen's \emph{Das
Prinzip der Infinitesimalmethode}.  While in \emph{Leibniz' System},
Cassirer concentrated on the historical figures of Leibniz and
Descartes, a few years later, he widened his perspective.  In his
monumental \emph{Das Erkenntnisproblem in der Neuzeit (The Problem of
Knowledge) } (Cassirer (1906--1950)) he became engaged in the huge
project of writing a comprehensive intellectual history of ideas
(\emph{Ideengeschichte}) of Western thought of the modern period that
he pursued during his entire lifetime in various forms: The last
volume of \emph{Das Erkenntnisproblem}, which eventually comprised
four bulky tomes, was published only postumously in 1950.  

Although in \emph{Leibniz' System} young Cassirer still followed
Cohen's lead in emphasizing the crucial importance of the
infinitesimal for modern science and mathematics, the reader may
notice an inclination toward relativizing its central role.  This
tendency gained momentum in later works such as \emph{Das
Erkenntnisproblem} (Cassirer 1906 - 1950), \emph{Substance and
Function} (Cassirer 1910 \cite{Cas10}), and \emph{The Philosophy of
Symbolic Forms} (Cassirer 1923 - 29).  But already in \emph{Leibniz'
System} we find the sweeping thesis that the central concept of modern
science is a concept of function%
\footnote{Here Cassirer is not referring to the notion of
function in a narrow mathematical sense, but rather to the functional
(or relational) account of science according to which entities are
secondary and functions or relations primary (see
Subsection~\ref{rel}).}
although it is not made clear how this assertion fits with the alleged
primacy of the concept of the infinitesimal. These tensions became
more evident in \emph{Substance and Function}.  Cassirer's failure to
toe the party-line on the primacy of the concept of the infinitesimal
over the concept of function was explicitly noted by Cohen; see
Section~\ref{five}.

For Cassirer the really modern character of Leibniz' thought was
encapsulated in the thesis that ``the real is conditioned by the
ideal.'' According to the Marburg idealism, this thesis was the key
that opened the possibility of a truly modern philosophy of science,
mathematics, and logic:
 
\begin{quotation}
From this perspective we can really understand Leibniz' tendency to
equate logic and mathematics in its true significance.  This equation
does not aim to constrict the rich content of mathematics in the form
of traditional logic.  Rather, it intends to bring about fundamental
reformation for logic.  Instead of being a theory of ``thought
forms'', logic is to become a science of objective knowledge
(\emph{gegenst\"andliche Erkenntnis}).  This transformation is
essentially due to its relation to mathematics: Mathematics turns out
to be the necessary mediation between the ideal logical principles and
the reality of nature. (Cassirer 1902 \cite[p.~113]{Cas02})
\end{quotation}

For Cassirer and Cohen, Leibniz was the one who opened the gate for an
idealist ``transcendental logic'' that later in the hands of Kant and
the Marburg neo-Kantians was to become a ``transcendental logic of
objective knowledge''.

\subsection{The Continuity principle} 
  
The essential means for overcoming the traditional Aristotelian
conception of logic as a theory of abstract ``thought forms'' toward a
contentful theory of objective knowledge was said to be Leibniz'
famous continuity principle.  The continuity principle was considered
as the most important device to unfold the general thesis that ``the
real is conditioned by the ideal'', mentioned in Subsection~\ref{ca}.
Consequently, the bulk of \emph{Leibniz' System} was dedicated to the
task of explicating this principle.

It should be kept in mind that for Cassirer and the Marburg
philosophers in general, continuity was not a property that some
things (or processes) have and others do not.  Cassirer described its
significance in the following terms:%
\footnote{Note that here Cassirer is using the concept of
``continuity'' as an ``operative concept'', as explained in
Subsection~\ref{rel}.}
\begin{quotation}
For the philosophy before Leibniz, continuity was essentially nothing
but the property of a thing or an attribute of a ready-made concept.
When it was understood in this way, one could attempt to refute or to
prove the claim that a certain thing or concept possessed it or lacked
it.
%
%\footnote{\label{f21}My predicament is that I don't understand this
%use of the term "predicament" here.  Perhaps you mean "predication" or
%something?}
%
This holds true of the \emph{synech\'es}
%
%$\sigma\upsilon\nu\epsilon\chi\epsilon\varsigma$ 
%
of the Eleates till Descartes' concept of continuous space.  Leibniz
overcomes this stance.  For him, the problem of continuity dissolves
in the problem of ``continuation''.  Continuity is no longer a
characteristic of a thing, but rather that of a development; not of a
concept, but of a method. (ibid., 153).
\end{quotation}

Therefore, it would perhaps have been more appropriate to refer to
this principle as the \emph{principle of continuation} rather than the
\emph{principle of continuity}.  In order to be understood as a
general conceptual achievement, \emph{continuation} needs to be
elaborated in the framework of a scientific methodology
(\emph{wissenschaftliches Verfahren}).  At this point the
infinitesimal and related concepts enter the stage.  As Cassirer was
eager to point out, the method of continuation first obtained its
deeper ``scientific'' meaning in the domain of geometry,
\begin{quotation}
where it designates the transition from point to line, from the line
to the area and so on.%
\footnote{From a purely mathematical viewpoint, Cassirer's formulation
is a bit unfortunate here because it sounds as if one is still dealing
with indivisibles rather than infinitesimals.  The difference between
them is that indivisibles were thought of as codimenion-1 entities
whereas infinitesimals were of the same dimension as the figure
composed of them.  This was the content of the major advance as
accomplished by Roberval, Torricelli, Wallis, Leibniz, and others, as
compared to earlier work by Archimedes and Cavalieri.}
Similarly, in mechanics and dynamics the method of continuation
describes the relation of the material point to the structures of
higher dimension.

In its true scientific generality the relation between an element and
the structure that results from its continuation corresponds to the
relation of a differential and its integral. In other words, the
``continuation'' is the methodical expression of the integration as a
continuous summation of infinitesimal moments. (ibid., 153 - 154)
\end{quotation}

To a modern reader, Cassirer's ``logical'' thesis may sound utterly
``metaphysical''.  Cassirer himself recognized that this general
characterization of the relation between the infinitesimal and the
real was in need of further clarification and scrutiny. For this task
he proposed to have a closer look at Leibniz' foundations of the
infinitesimal calculus.  According to Cassirer, the key to
understanding the true novelty of Leibniz' account resided in the
Leibnizian concept of motion. The crucial point was not to view motion
as an empirical concept stemming from the realm of empirical
experience.  For Leibniz, motion was always \emph{continuous} motion,
i.e., the expression of a unifying principle of conceptual
construction. Leibnizian motion was not something empirically given,
but something conceptually constructed (cf. ibid., 156). Invoking the
idealist principle that ``the real is clarified by the ideal'', still
another way of expressing this may be that, the concept of continuity,
in Leibniz's view, belongs to the realm of the ideal.
 
Probably the most famous expression of this view can be found in
Leibniz's letter to Varignon from which Cassirer quoted in his
\emph{Leibniz' System} (ibid., p. 188/189) and elsewhere:

\begin{quote}
\ldots one can say in general that, though continuity is something
ideal and there is never anything in nature with perfectly uniform
parts, the real in turn, never ceases to be governed perfectly by the
ideal and the abstract \ldots (Leibniz 1702 \cite{Le02})
\end{quote}

Later, in the same letter, Leibniz explicitly stated that not only
continuity but also infinitesimals have the capacity to ``govern the
real perfectly":

\begin{quote}
So it can be said that infinites and infinitesimals are grounded in
such a way that everything in geometry, even in nature, takes place af
if they were perfect realities. (ibid.) 
\end{quote}

On the other hand, in the very same letter he asserted: 

\begin{quote}
I'm not myself persuaded that it is necessary to consider our
infinities and infinitesimals as something other than ideal things
(choses ideales) or wellfounded fictions (fictions bien
fond\'ees). (ibid.)  
\end{quote}

It is far from clear, however, how precisely the relation between
``ideal things" and ``(wellfounded) fictions" is to be thought and how
infinitesimals as ``fictions" could have this power of governing
perfectly.  This difficulty was already observed by Cassirer, who even
contended that Leibniz's notion of a ``fiction bien fond\'ee" had an
air of paradox (Cassirer 1902 \cite[p.~187ff]{Cas02}).  For a
contemporary survey of the debate concerning this issue of Leibnizian
scholarship the reader may consult (Sherry and Katz 2013 \cite{SK}).

%\subsection{A random break: Cassirer on Leibniz, continuity, and function}

Cassirer himself contended that the doctrines of the Marburg school,
in particular his account of the role of idealizations in science and
mathematics, might help to overcome the remaining obscurities that
still beset Leibniz's account. According to him, for this endeavor it
was essential to properly understand Leibniz' concept of motion that
underlied his dynamical conception of geometry. From the viewpoint of
classical Euclidean geometry, one may suspect that introducing the
concept of motion into geometry amounts to an illicit confusion of
pure geometry and empirical science. Cassirer vigorously argued
against this interpretation:

\begin{quotation}
It is not a systematical infringement
%
%($\mu\epsilon\tau\alpha\beta\alpha\sigma\iota\varsigma$
%$\epsilon\iota\varsigma$ $\alpha\lambda\lambda$o
%$\gamma\epsilon\nu$o$\varsigma$) 
%
to integrate the concept of motion into geometry.  The concept that is
dealt with here, is not from physics, but from logic: It denotes the
conceptual continuation of the ``principle'' that was expressed in the
concept of continuation.  Thereby the concept of motion is separated
from its empirical context and allocated in the area of pure and
eternal ``forms''.  \ldots The general achievement of the concept of
motion resides in the formulation of the thought that the extensional
being has to be constituted from an original lawful determination that
preceeds it as its logical \emph{prius}.%
\footnote{A \emph{prius} is something that comes before or preceeds
something else in some respect. The term was used in the philosophical
jargon of the 19th century. Today it seems to be an outmoded
Latinism.}
(Cassirer 1902 \cite[pp.~156-157]{Cas02})
\end{quotation}

In other words, motion in the realm of science is always lawful
motion.  The purely conceptual character of the concept of motion, its
non-extensionality is shown by the concept of the differential
(ibid. 157). This is not to say that Cassirer was not aware of the
existence of non-continuous and perhaps even continuous but
non-differentiable functions which had been vigorously discussed among
mathematicians since Cauchy in the 1820s (cf.~Hankel 1882
\cite{Han82}).  Indeed, he explicitly pointed out that a more general
notion of a function \`a la Dirichlet, perfectly made sense from a
purely mathematical point of view.%
\footnote{Compare footnote 8.}
He only objected that such general functions were not meaningful for
the determination of real processes of nature (cf. 217).  He went on
to declare that the meaningfulness of the limiting processes of the
calculus ``demonstrated'' that the concept of (lawful) motion in
nature was a non-empirical, logical notion:

\begin{quotation}
If the transition to the quantitative zero does not eliminate the
lawful character of the magnitude this is evidence that it (i.e.,
lawfulness) is not grounded in a quantitative principle. The magnitude
must first disappear from our sensual perception before we can
recognize its determinateness in the pure concept (Cassirer 1902,
157).
\end{quotation}
Here Cassirer does not distinguish between continuity and
differentiability.  As we shall see in a moment, this conflation
enabled him to combine Cohen's ``infinitesimal-centered'' account with
his own ``function-oriented" one in an elegant but somewhat dubious
way.

For the contemporary reader this passage may sound opaque, to put it
mildly. We propose the following interpretation. Scientifically
meaningful magnitudes obey continuous motions, i.e., motions the law
of which could be described by a differentiable function $f$.  In
calculating the derivative of~$f$, expressions such
as~$\lim\frac{f(x)-f(x')}{x-x'}$ occur.  These contain ``quantitative
zeroes'' (if $x'$ approximated $x$).  This means that such expressions
have no direct quantitative meaning.  In particular they could not be
perceived or experienced in any reasonable way.  Nevertheless,
conceptually, the calculation of the derivative makes perfect sense.
Hence, the magnitude could be recognized as a meaningful and
determinate magnitude, only after it had been submitted to a
conceptual process (derivation) that eliminated all its sensual
qualities. Continuity was thereby asserted to be a necessary
presupposition for the constitution of nature as a possible object of
rational investigation:
   
\begin{quotation}
Continuity is a necessary presupposition for the existence of a
mutually 1-1 relation between two series of change
(Ver\"anderungsreihen). This strand of thought is first formulated in
Leibniz' best known formulation of the principle of continuity:
``Datis ordinatis etiam quaesita sunt ordinata''.  The ``data'' denote
the hypothetical conditions from which we start; the ``quaesita'' are
the series of the conditioned that we look for. The order is thought
as a law that determines the transition inside the two series in a
continuous fashion. (ibid., 211 - 212)
\end{quotation} 

Cassirer went on to contend that the standard ``epsilon-delta''
definition of continuity was merely a mild reformulation of Leibniz'
original characterization of continuity ``Datis ordinatis etiam
quaesita sunt ordinata'' (cf. ibid., 215).  Thereby he could conclude
that there was an intimate relation between Leibniz' continuity
principle, the modern epsilon-delta definition of continuity, and
Cohen's infinitesimal-centered account.  Modern mathematics has shown
that the relation is more complicated than Cassirer might have
thought.

Nevertheless, despite its allegedly close relation to the modern
mathematical concepts, according to Cassirer, Leibniz' continuity
principle should not be understood as a mere mathematical definition.
Here Cassirer is treading on somewhat dangerous ground.  Leibniz's law
of continuity had several meanings. It was not a single concept but
rather a family of concepts; see for example (Jorgensen 2009
\cite{Jor09}).  The problem is not merely the fact that the
mathematical concept of continuity is not general enough to encompass
Leibniz's concept.  Rather, the law of continuity invokes several
related concepts in a chain where the concepts at the extrema may be
unrelated to each other at all.  Perhaps one may say that Leibniz used
this concept as an `operative' one (see Subsection~\ref{rel}).

It would be a misunderstanding to read it simply as the claim that the
processes of nature should be viewed as ready-made entities that were
to be described in terms of continuous (or differentiable) functions:

\begin{quote}
The requirement of conceiving of nature ultimately as a complex of
continuous functions (\emph{Inbegriff stetiger Funktionen}) would not
make sense, if the task of knowledge would be to reproduce a
ready-made material descriptively.  Continuity obtains its meaning
only if it is conceived of as a basic act of the mind through which
the subject conditions the object. (ibid., 218)
\end{quote}
 
The principle of continuity should be understood as a guiding maxim
for the evolution of scientific concepts that urges us to seek ever
more profound systematic connections among 
them.  Concepts should be connected in a uniform conceptual system and
each concept should be transformable into each other continuously. For
the Marburg school Leibniz' principle of continuity was the
philosophical expression of one of the basic moments of modern
science. It asserted that the well-defined and determinate character
of scientific concepts did not reside in their isolation but in the
lawfulness of their transitions (cf.~(ibid. 201)).

Although Cassirer, faithfully following in the footsteps of his master
Cohen, hailed Leibniz as the genius who made explicit for the first
time the principle of continuity as a fundamental principle of modern
science, already in \emph{Leibniz' System} his assessment of Leibniz's
achievements went in directions other than that of Cohen.  The
emerging differences between Cohen and Cassirer concerned the relation
between the concepts of the infinitesimal and function.  While Cohen
emphasized that the ``infinitesimal calculus (of Leibniz) had placed
the concept of function, conceived of as a law of interdependency
between two variable magnitudes, in the center of the methodology of
mathematics" (Cohen 1902 \cite[p. 239]{Co02}), Cassirer put less
emphasis on the role of the infinitesimal as a conceptual base for the
concept of function. He generally praised Leibniz as an early partisan
of a ``functional" or ``relational" worldview without mentioning
infinitesimals at all:

\begin{quotation}
If one understands by ``substantial" worldview the conception
according to which all beings and occurrences can be traced back to
ultimate, rigid, absolute ``things'', then Leibniz' philosophy is
strictly opposed to this standpoint.  The tendency of Leibniz'
philosophy that from now on will prevail in the ongoing progress of
idealism points at a replacement of the older concept of \emph{being}
by the concept of \emph{function}. (ibid., 486) [emphasis added--the
authors]
\end{quotation}

In the evolution of Cassirer's own thought, this functional
interpretation of the principle of continuity gained ever greater
momentum and superseded the infinitesimal interpretation eventually
leading to certain discrepancies with Cohen that surfaced in Cohen's
letter to Cassirer dating august 24, 1910 (see Section~\ref{five}).

\subsection{Gawronsky's \emph{The Judgment of Reality}}

The second sustained effort to make sense of Cohen's approach is due
to Dimitry Gawronsky (1883 - 1955).  His dissertation under Cohen and
Natorp was entitled \emph{Das Urteil der Realit\"at} (Gawronsky 1910
\cite{Gaw10}).  Although he was a close friend of Cassirer's, in the
emerging discrepancies between Cohen and Cassirer on the relation
between Cohen's `infinitesimal' and Cassirer's `functional' approach,
he sought to find a mediating position between the two but eventually
sided rather with Cohen than Cassirer.  For this issue, two works of
Gawronsky are relevant.  Besides his already mentioned dissertation
\emph{Das Urteil der Realit\"at} (henceforth \emph{Urteil}), we have
also his contribution \emph{Das Kontinuit\"atsprinzip bei Poncelet}
(Gawronsky 1912 \cite{Gaw12}) to a \emph{Festschrift} dedicated to
Cohen on the occasion of his 70th anniversary in 1912.  Today,
Gawronsky's philosophical work has fallen into almost complete
oblivion.  Yet he was an important figure in the internal debate that
took place within the Marburg school on matters infinitesimal between
Cassirer, Natorp, and Cohen in the early years of the 20th century.

Gawronsky took upon himself the difficult task of updating Cohen's
infinitesimal account, defending it against the less than orthodox
accounts of Cassirer and Natorp. Probably his best known work among
Cassirer scholars is the biographical article \emph{Ernst Cassirer:
His Life and His Work} (Gawronsky 1949 \cite{Gaw49}) that appeared as
a contribution to the Schilpp volume dedicated to Cassirer.  As far as
we know, the only contemporary discussion of Gawronsky's work and his
role as a vigorous (although not uncritical) defender of Cohen's
position against Cassirer (and, to a lesser extent, Natorp) is Massimo
Ferrari's paper \emph{Dimitry Gawronsky and Ernst Cassirer: On the
History of the Marburg School between Germany and Russia} (Ferrari
2010 \cite{Fe10}) published in Russian.

In contrast to Cohen, Gawronsky was fully competent in matters of
contemporary mathematics.  He discussed the achievements of Bolzano,
Grassmann, Cantor, Weierstrass, Veronese, and Dedekind with evident
expertise. Moreover, Gawronsky expressed a positive appreciation of
the limit method (cf. Ferrari [35, p. 249]).

Nevertheless, in contrast to Cassirer, Gawronsky sought to leave the
philosophical core of Cohen's `infinitesimal-centered' account intact.
Hence, with respect to the infinitesimal approach he, rather than
Cassirer or Natorp, may be considered as Cohen's true heir.  His
dissertation \emph{Das Urteil der Realit\"at und seine mathematischen
Voraussetzungen}, literally `The judgment of reality and its
mathematical premises', may be regarded as the only serious attempt of
amending Cohen's rather obscure pseudo-Kantian table of judgments.

In order to arrive at a better understanding of Gawronsky's
\emph{Urteil}, we have to recall that, as Gawronsky explains at the
end of \emph{Urteil}, the `judgment of reality' refers to a
distinction already made by Kant (and later modifed by Cohen) that can
only be translated with difficulty into English.  This is the
distinction between \emph{Wirklichkeit} and \emph{Realit\"at} that are
both usually translated as \emph{reality}.  Roughly, \emph{reality} in
the sense of \emph{Realit\"at} is to mean `the systematic knowledge of
nature as it arises from the chaos of the immediately sensed'
(\emph{unmittelbares Empfinden}) (Gawronsky 1910
\cite[p.~107]{Gaw10}).  Then the main thesis of \emph{Urteil} is that
the infinitesimal calculus plays a crucial role in the systematic
knowledge of nature (or empirical reality) insofar as
\begin{quotation}
the basic problem of objective empirical knowledge is the problem of
change.  But we only understand change if we obtain complete knowledge
of the law that generates change, i.e., only if we can pursue the
effect of the generating law in every infinitely small element of this
change$\ldots$ And exactly this is achieved by the infinitesimal
analysis. (Gawronsky 1910, \cite[p.~105]{Gaw10}).
\end{quotation}

In line with Cohen, Gawronsky asserted that ``there is no other way to
formulate and to justify the laws of nature than the infinitely
small.'' (ibid.)

By identifying ``reality'' with the systematized knowledge of nature,
Gawronsky saw ``reality'' as a ``logical method'' whose essence ``was
the assumption of the existence of a generating law''.  This ``logical
method'' came along in two different ways, namely ``by the method of
number and by the method of the infinitesimal'' (108).  Both methods
are carried out in the same three steps:

\begin{enumerate}
\item
Positing (\emph{Setzung})
\item
Infinite repetition (\emph{unendliche Wiederholung})
\item
Actual synthesis in a higher unity (\emph{aktuale Zusammenfassung in
einer h\"oheren Allheit})
\end{enumerate}

Gawronsky's attempt to construe an analogy between the two methods is
apparently based on the idea, which he shared with Dedekind and
Cassirer, that the essence of numbers resided in their ordinal
structure.  More precisely, according to Gawronsky the conceptual
generation of the natural numbers proceeded by first positing the unit
`1' and then applying the generating principle of the successor
function, thereby constituting the other natural numbers.  This
construction, however, had to be `completed' by an ``actual synthesis
resulting in a higher unity".  Or, formulated negatively, Gawronsky
was not content with simply asserting that this repetition could go on
and on leading to ever larger natural numbers.  Rather, one had to
look for a higher synthesis.

This was achieved, Gawronsky contended, by Cantor's theory of infinite
ordinals.  More precisely, Gawronsky conceived of Cantor's positing of
the first infinite ordinal~$\omega$ as the sought-for completion or
synthesis. This completion of the natural numbers in terms of the
first infinite ordinal~$\omega$, however, was not simply the end of
the constitution process of pure thought.  On the contrary, it was
just the beginning of a new stage in that it gave rise to a new
infinite series generated by a new generating principle~$\omega, \,
\omega +1, \, \omega+2, \ldots$.  This new series, then, had to be
conceptually completed by positing 2$\omega$, which served as the
starting point for a new series 2$\omega, \, 2\omega +1, \, 2\omega+2,
\ldots$.and so on. The determination of the limit of an infinite
arithmetical series such as $3.1, 3.14, 3.141, 3.1415, \ldots$
converging to~$\pi$ had a similar conceptual structure, and even for
the calculation of differentials and derivations Gawronsky assumed an
analogous conceptual structure.  According to him, they all followed
the three-tiered pattern of `positing', `repetition', and `actual
synthesis'.

A modern mathematician may view Gawronsky's contention merely as the
recognition that the construction of both the real numbers and
infinitesimals involves infinitary constructions, a point made much
later also by Robinson.  But Gawronsky, following his master Cohen,
made much more of it.  According to him, the usage of infinitary
constructions revealed the very essence of both (empirical) science
and mathematics as being based on ``pure thought''.  The
transcendental analysis of science revealed that the origins of these
methods were to be found in the `judgment of reality'.  Something that
could not be counted%
\footnote{\label{f31}It should be noticed that Gawronsky here relied
on a rather broad concept of counting that not only included
``ordinary" counting but also various kinds of ``infinite completion
of counting''.}
or differentiated, was not ``real" in the sense that it could not
possibly be the object of scientific knowledge.%
\footnote{This has the apparently paradoxical consequence that there
may be something real - in the sense of \emph{wirklich} - that is not
\emph{real}.  This apparent contradiction does not threaten in the
original German and is avoided if one carefully distinguishes between
the two meanings of ``real" in Kant's language.}
Differentiating and counting were the two basic methods of scientific
knowledge. From a modern point of view, this may be a somewhat narrow
and outdated characterization of the conceptual apparatus used in
science and mathematics, but it certainly makes sense.%
%
%\footnote{\label{f27}I suggest deleting this material: But Gawronsky's
%contention went further.  He claimed to have succeeded, guided by
%Cohen's insights in \emph{Prinzip} and \emph{Logik}, in determining
%the ``logical judgment" underlying these two methods, namely, ``the
%judgment of reality, whose essence we described as that of the
%generating law.  Thereby the transcendental deduction of the judgment
%of reality had been achieved, i.e., it had been shown that it is a
%basic method of scientific knowledge and at the same time a necessary
%component of the system of transcendental logic as the logic of pure
%knowledge'' (Gawronsky 1910 \cite[p.~109]{Gaw10}).  Even a sympathetic
%reader might complain that this passage is less than crystal-clear.}
%

Let us now examine \emph{Das Kontinuit\"atsprinzip bei Poncelet}
(Gawronsky 1912).  The main aim of this work was to elucidate Cohen's
dictum that \emph{continuity} is a basic law of thought
(\emph{Denkgesetz}) (Cohen 1902 \cite[p.~76]{Co02}).  In \emph{Logik},
Cohen traced the principle of continuity back to Leibniz.
Furthermore, he offered the following high-sounding explication of the
role of continuity in the ongoing process of philosophical and
scientific thought:
\begin{quotation}
Continuity is a law of thought.  It is the law of thought of the
connection which enables the generation of the unity of knowledge and
thereby the unity of the object of knowledge. É Continuity as a law
of thought garantees the connection of all methods and disciplines of
mathematized empirical science (\emph{mathematische
Naturwissenschaft}).  This law is therefore of crucial importance for
the thinking of knowledge.  Continuity is the law of knowledge.
Continuity characterizes the basic feature of thought (Cohen 1902, 76,
77).
\end{quotation}
Gawronsky's aim in \emph{Das Kontinuit\"atsprinzip} was to confirm and
elucidate Cohen's global thesis on the central role of the principle
of continuity for the evolution of scientific thought, by studying its
role in the development of 19th century geometry exemplified in the
work mainly of French geometers like Poncelet, Chasles, Carnot, and
others.  In line with Cohen, Gawronsky points out that this principle
is not an achievement of 19th century science but a basic feature of
all scientific thought. What is new, according to Gawronsky, is the
way Poncelet applied the principle.  The crucial point is not that new
objects are subsumed under the known theorems and relations
(cf. Gawronsky 1912, 69) but that entire systems of theorems and
relations themselves are `continuously' modified and generalized
(ibid. 71).  This new interpretation of the principle entails that it
must not be understood as an argument that generates mathematically
secure results but rather as a heuristic principle that helps one find
novel concepts whose relevant connections have yet to be secured by
other means:
\begin{quotation}
We see that the formation of concepts that is determined and guided by
the principle of continuity cannot be completely justified by that
principle alone. Rather, a subsequent check has to be carried out in
order to determine the value of every newly introduced
concept. (Gawronsky 1912 \cite[p.~73]{Gaw12})
\end{quotation}
For Gawronsky, the anticipatory and heuristic character of the
principle of continuity in mathematics, as it was employed by Poncelet
and others, was essential for discovering its true logical base that
comes to the fore when we compare it with its ``prototype'' as it
appears in the realm of pure thought:
\begin{quotation}
Since the discovery of the principle of continuity, Hermann Cohen in
his \emph{Logik der reinen Erkenntnis} was the first who sought not
only to evaluate the achievements of this principle in a comprehensive
way, but also to determine its systematic significance, to introduce
it in the system of pure thought and to render precise its position.
He was the one who recognized this priniciple as a rather general and
basic method of scientific thought, identifying it as a law of thought
of knowledge. (Gawronsky 1912, 74)
\end{quotation}
From the viewpoint of Cohen's neo-Kantian approach, it is then the
task of philosophy of science to integrate the issue of a purely
mathematical evaluation of the principle of continuity, as it was
understood by Poncelet and his contemporaries, into the general agenda
of the transcendental logic of science that treats the basic methods
of scientific knowledge \emph{\"Uberhaupt} (ibid. 76).  For this
purpose the philosopher must not rely on idle metaphysical
speculations but has to know how this principle is actually applied in
scientific practice.  As an example of how this may be achieved in the
case of geometry, Gawronsky discusses in detail two principles that
Poncelet introduced in modern geometry, namely, the principles of
central projections and his ``th\'eorie des polaires reciproques"
(ibid., 76ff). This leads him to the conclusion that also Klein's
\emph{Erlanger Programm}, which proposes to define the essential
properties of geometrical objects as invariants of certain
transformations groups, can be unterstood as a realization of the
principle of continuity (78).  In sum, one may contend that
Gawronsky's \emph{Das Kontinuit\"atsprinzip bei Poncelet} offers a
knowledgeable and not implausible narrative of the development and
significance of the principle of continuity in 19th century geometry.
This can be taken as indirect evidence that Husserl overstated his
case when he summarily dismissed Cohen's account of the principle of
continuity as a basic law of thought as ``profound nonsense''.
Admittedly, it often takes considerable effort to distill some meaning
out of Cohen's obscure prose, but the attempt to rescue at least some
parts of Cohen's transcendental logic as presented in \emph{Die Logik
der reinen Erkenntniss} cannot be bluntly dismissed in Husserl's
fashion.

\subsection
{Natorp's \emph{The Logical Foundations of the Exact Sciences}}

In his \emph{Die logischen Grundlagen der exakten Wissenschaften}
(Natorp 1910 \cite{Na10}) the second leader of the Marburg school,
Cohen's friend and colleague Paul Natorp also sought to come to terms
with the problem of infinitesimals.  He sought to develop a
philosophically founded synthesis of two antagonistic, or at least
very different, mathematical programs for the foundations of analysis.
These are, on the one hand, the CDW program, which in Natorp's day had
nearly achieved the status of a ruling orthodoxy, and on the other,
the maverick, ``infinitesimal-friendly'' program of the Italian
mathematician Giuseppe Veronese, put forward in his \emph{Fondamenti
di geometria} (Veronese 1892).%
\footnote{See P. Cant\`u \cite{Cantu} for an extensive bibliography on
Veronese.}

Natorp's contribution to the Marburg neo-Kantian philosophy in
general, and to mathematics in particular, has been usually neglected
compared to the better known works of Cohen and Cassirer.  This may be
considered, historically speaking, as an injustice, in particular with
regard to the issue of infinitesimals.  In Natorp's \emph{Die
logischen Grundlagen der exakten Wissenschaften} (Natorp 1910) we find
the most elaborate and most complete discussion of limits and
infinitesimals that any neo-Kantian philosopher ever published.%
\footnote{\emph{Die Logischen Grundlagen} has never been translated
into English.  Even in German-speaking philosophy Natorp always
remained in the shadow of the more brilliant Cassirer.  
%
%We know of only one recent interpretation of Natorp's account of
%infinitesimals, limits, and continuity, namely Peiffer-Reuter's
%articles \cite{Pei89}, \cite{Pei92}.%
}

Indeed, Natorp dedicated two chapters (namely, chapter~III, 98 - 159,
and chapter~IV, 160 - 224) of \emph{Logische Grundlagen} to a detailed
criticism of the accounts of the various concepts of number as put
forward by Frege, Dedekind, Cantor, Weierstrass, Pasch, and Veronese.
Furthermore, Natorp believed that the accounts of Cantor and Veronese
are compatible, and viewed Veronese as the ``most eminent successor of
Cantor'' \cite[p.~171]{Na10}.  He appears to have held that the
differences between them were only technical differences of no
conceptual and philosophical relevance.  This was certainly an error,
as Cantor and Veronese were well aware of the fact that their accounts
differed in essential ways.  Cantor rejected Veronese's numbers.
Veronese was more tolerant, seeing Cantor's transfinite numbers and
Veronese's own transarchimedean%
\footnote{This somewhat unusual coinage was utilized by
Laugwitz~\cite[p.~104]{Lau02} and
Peiffer-Reuter~\cite[p.~124f.]{Pei89}.}
infinitely large numbers as two admissible, but nevertheless quite
different types of mathematical entities.

An analysis of Veronese's account would go beyond the scope of the
present text.  What we wish to explore is Natorp's philosophical
motivation that inspired him to engage in the risky endeavor of
sketching an all-embracing panoramic view of the landscape of the
various kinds of numbers%
\footnote{\label{ff2}In retrospect, one can assert that Natorp's
all-embracing vision was not merely a philosopher's pipe dream.  It
gained some mathematical substance over the decades.  Thus, Kanovei
and Reeken proved that there is a certain (class-size) structure
${}^*\!V$ which is $\kappa^+$-saturated for every cardinal $\kappa$,
together with an elementary embedding of the ZFC set universe $V$ into
${}^*\!V$ (Kanovei and Reeken 2004 \cite[Theorem~4.3.17, p.~151]{kr}).
(Ehrlich 2012 \cite{Eh11}), working in Von Neumann--Bernays--G\"odel
set theory with Global Choice, showed that a maximal (class) hyperreal
field is isomorphic to the maximal \emph{surreal} field.}
and their calculi (cf. Natorp~\cite{Na10}).

In line with all neo-Kantians, Natorp contended that intuition alone
could not provide a foundation for our knowledge about infinitesimals,
limits, and continuity.  As he unequivocally put it: ``Intuition
cannot be a foundation for continuity, neither for space nor for time"
\cite[chapter~IV, \S4, p.~175]{Na10}.

According to him, the only relevant factor was ``pure thought'',
i.e. the ``transcendental method'' or ``transcendental logic'' (see
Subsection~\ref{31}).  Indeed, Natorp's insistence on the crucial role
of an unfettered investigation of infinitary objects in applying the
``transcendental method'' gave his interpretation of infinitesimals in
particular and numbers in general its specific flavor:%
%
%commented out material left for historical reasons, to be ignored by
%both authors: The quote provided does not really bear out the
%infinitary promise of the previous sentence.  There is nothing here
%explicitly about infinity.
%
\begin{quote}
Numbers must not have any other basis than the laws of pure thought.
(Natorp \cite{Na10}, ch. IV. \S4, p.176).
\end{quote}
His insistence on the infinitary character of the never-ending road of
the ``transcendental method'' led him to criticize attempts to base
the concept of number on something ``finite'' that only at a later
stage of the conceptual evolution was overcome in favor of something
``infinite''.%
\footnote{Natorp's insistence on the thoroughly infinite character of
numbers, and the resulting emphasis on infinitesimals, was noted not
only by his fellow philosophers but also by philosophically inclined
mathematicians, regardless of whether they belonged to the CDW camp or
were sympathetic to infinitesimals; see Fraenkel \cite[footnote on
pp.~50-51]{Fra12} and Robinson \cite[p.~278]{Ro66}.}
For him, from the very start, the concept of number was soaked with
the infinite.  Hence, for him, any attempt to conceive of the
``finite'' rational numbers as a more solid base for allegedly the
more elusive real numbers was philosophically mistaken.  He appears
even to have blamed Dedekind for having succumbed to such a temptation
to some extent.

Natorp also came to formulate a perceptive criticism of a reductive
aspect of Dedekind's approach of introducing new numbers as ``cuts''
of the set of rational numbers (cf. Dedekind 1872).  Natorp pointed
out that Dedekind assumed without justification that every cut%
\footnote{We will ignore the technical issue of cuts defined by the
rationals themselves.}
corresponds to exactly one non-rational number, which according to
Natorp was a \emph{petitio principii}.  Of course, two distinct
non-rational numbers which correspond to one and the same cut, cannot
differ by a finite number.  However, they could still differ by an
infinitesimal.%
\footnote{This is indeed a mathematically coherent possibility as is
shown already by the Levi-Civita fields~\cite{Le} developed at about
the same time.  In fact, any proper ordered field extension of the
reals will have this property (for example, the hyperreal numbers).}
Thus, Dedekind's axiom that every cut corresponds to a \emph{unique}
number is an assumption that can be challenged.

Indeed, Natorp was the only neo-Kantian philosopher who ever
explicitly evoked a suspension of the Archimedean axiom as an
essential ingredient for any infinitesimal account
(cf. \cite[p.~169ff]{Na10}).  He evoked the Archimedean principle in
his criticism of Cantor's alleged ``proof'' of the non-existence of
infinitesimals (ibid.).%
\footnote{Natorp claimed that Cantor's ``proof'' was flawed for rather
trivial reasons.  His remarks are, however, too sketchy to be properly
evaluated.  Cantor was a competent mathematician and although he did
make some mistakes, they were rather subtle ones.  On the other hand,
Natorp's competence in matters mathematical was that of an educated
layman.  Cantor's errors are analyzed in detail by Ehrlich~\cite{Eh06}
and Moore~\cite{Moo}; see also (Proietti 2008 \cite{Pr}).}

Natorp further points out that already the classical founding fathers
of the concept of the infinitesimal, namely, Leibniz and Newton were
well aware of the special status of the infinitesimal, often referred
to by means of the modifier ``intensive''.  The latter term seems to
have had a meaning close to ``non-Archimedean'' to Natorp:
\begin{quote}
Already Galileo speaks of \em{infinita non quanta}\em; Leibniz
contends the infinitesimal as \em{praeter extensionem, imo extensione
prius}\em, for Newton the infinitesimal ``moments'' are not
\em{quantitates finitae}\em, but \em{principia iamiam nascentia
finitarum magnitudinum}\em; and Kant explains the infinitesimal
through the intensive magnitude that include the base (\em{Grund}\em)
for the extensive magnitudes, but is itself not extensive
(\cite[p.~170]{Na10}).%
 
\end{quote}

Natorp further criticized Dedekind's assumption that the totality of
all cuts is `given' somehow although the converging rational series
are to be considered just as procedures that allow us to approximate
the irrational limit numbers \emph{if} they exist.

Natorp's criticism of the idea that the set of Dedekind cuts is
``given" could appear to be related to a criticism often voiced at the
time, and that can be expressed in modern mathematical terms as
follows.  To say that the set of cuts is ``given'' is to make certain
foundational assumptions, such as the axiom of infinity, usually
accompanied by the classical interpretation of the existence
quantifier, typically involving the law of excluded middle.  If so,
Natorp's criticism of Dedekind seems to echo criticisms of
mathematicians like Kronecker and Brouwer.  This, however, would be a
misinterpretation.  Natorp was neither an intuitionist \`a la Brouwer
nor a platonist realist.  Rather, Natorp was a \emph{critical}
idealist.  For him, the ideal was neither something given ``out
there'' nor did he require that it could be intuited in some way or
other. For him, as well as for Cassirer, the existence of the ideal
resided only in its function.  This led him to perceive similarities
between Dedekind's approach and that of du Bois-Reymond:

\begin{quotation}
Apparently this kind of argumentation [i.e., Dedekind's---the authors]
is based on a way of thought that is manifest in P. du Bois-Reymond's
\emph{Allgemeine Funktionentheorie} (\emph{General Theory of
Functions}).  This author introduces infinity and continuity by
nothing short of an assumption (which he himself calls ``idealist''
but which is actually ``realistic'' in the sense of medieval
scholasticism) that can often be seen in arithmeticians: namely, that
the objects of mathematics exist in-themselves, and that these objects
may have properties which our - always finite - human thinking cannot
fully grasp. (Natorp 1910, 180)
\end{quotation} 

Natorp insisted on the thesis that the existence of mathematical
objects can reasonably only mean that they are based on the
mathematical thought (ibid.). Something that escapes mathematical
thought, does not exist mathematically. In other words, du
Bois-Reymond's conception of the continuum as something ``that cannot
be thought" did not make sense:%
\footnote{At about the same time, a similar, but more elaborate
critique of du Bois-Reymond's ``empirist-idealist dialectics'' was put
forward in (Cassirer 1910 \cite[p.~122ff]{Cas10}).}
\begin{quotation}
What cannot be justified by mathematical thought, must not be posited
by mathematics" (ibid. 180).
\end{quotation}

Despite certain alleged shortcomings in Dedekind's cut approach,
Natorp saw Dedekind as being on the right track.  The merit of having
revealed the true kernel of Dedekind's method, is ascribed by Natorp
to Weierstrass, Cantor, Pasch, and Veronese.  He sees the basic flaw
in Dedekind in the fact that Dedekind started his construction with
the ``finite'', fully understood rational numbers, whereas the
irrational numbers were considered as something derived.%
\footnote{This criticism of Natorp's is related to Cassirer's
insistence on the ``ontological equality'' for the new entities (in
this case, numbers) being introduced.}
According to him, Weierstrass and Cantor made the decisive conceptual
step of taking the infinite convergent ``series" itself (Cantor's
\emph{Fundamentalreihe}) as a proper mathematical object to be
considered in its own right.  Then there is no longer any reason to
distinguish between a ``series'' and its limit (cf. Natorp 1910
\cite[p.~182]{Na10}).

\section{From infinitesimals to functional concepts}
\label{five}

Cassirer sought to imbed Cohen's infinitesimals in a larger framework
built upon his functional approach (cf. Cassirer 1912 \cite{Cas12}),
but faithfully followed the general doctrine of the Marburg school not
to dismiss infinitesimals as inconsistent pseudo-concepts. Unlike
Cohen, however, he was interested not so much in infinitesimals
\emph{per se} as in a more profound and more precise understanding of
the conceptual evolution as it took place in the modern sciences, in
particular in physics and mathematics.  For this purpose, he was led
to a substantial recasting of the neo-Kantian framework of philosophy
of science and mathematics as it had been designed by Cohen in
\emph{Prinzip} and \emph{Die Logik der reinen Erkenntnis}.

In a nutshell, this amounted to placing the concept of \emph{function}
or \emph{relation} center stage, rather than that of the
infinitesimal.  This project of Cassirer's began with his dissertation
\emph{Leibniz' System in seinen wissenschaftlichen Grundlagen}
(Cassirer 1902 \cite{Cas02}).  The project was succinctly presented in
\emph{Kant und die moderne Mathematik} (Cassirer 1907 \cite{Cas07}),
and eventually culminated in \emph{Substanzbegriff und
Funktionsbegriff} (Cassirer 1910) and \emph{Die Philosophie der
symbolischen Formen III} (Cassirer 1929 \cite{Cas23}).

Cohen was not entirely happy with this development as shown by his
letter to Cassirer dating from August 24 of 1910.  Here Cohen first
heaps lavish praise on his most brilliant disciple:
\begin{quotation} 
I heartily congratulate you and our entire community on your new and
great achievement [i.e. the publication of SF---the authors].  If I
shall not be able to write the second part of my \emph{Logik}, no harm
will be done to our common cause.%
\footnote{Cohen never published a second part of \emph{Logik}.}
\end{quotation}
Cohen then continues with the second half of his comment:
\begin{quotation}
Yet, after my first reading of your book I still cannot discard as
wrong what I told you in Marburg: you put the center of gravity upon
the concept of relation and you believe that you have accomplished
with the help of this concept the idealization of all materiality. You
let even slip the remark that [the concept of relation] is a
category. (\ldots) Yet it is a category only insofar as it is a
function, and function unavoidably demands the infinitesimal element
in which alone the root of the ideal reality can be found.
\end{quotation}

Among the recent interpretations of the relationship between Cohen's
infinitesimals and Cassirer's relational approach, one can find
conflicting views, both of which take this letter as their main piece
of evidence.  In his 2003 article \emph{Hermann Cohen's Das Prinzip
der Infinitesimalmethode, Ernst Cassirer, and the Politics of Science
in Wilhelmine Germany}, Moynahan \cite{Moy} put forward the thesis
that Cassirer's relational account in SF should be understood as a
more or less straightforward clarification of Cohen's \emph{Prinzip}.
As evidence for his claim, he quotes Cohen's 1910 letter to Cassirer.
However, Moynahan only quotes the first half.  The second half of the
letter, in which Cohen pointed out the profound differences between SF
and \emph{Prinzip}, is not reproduced by Moynahan; see
\cite[p.~40]{Moy}.

In a more subtle and indirect way than Moynahan, recently Seidengart
(2012) also argued that Cohen and Cassirer essentially agreed on the
primordial role of the concept of the infinitesimal for modern science
and its proper philosophical understanding.

To bring home his point Seidengart first reminded the reader that for
Cohen the concept of the infinitesimal had to occupy center stage in
any logic of modern science deserving of its name, since
``infinitesimal analysis was the legitimate device of the mathematical
science of nature."  (Cohen, 1902, p.~30).  According to Cohen,
Leibniz, as the inventor of the infinitesimal calculus, was the one
who brought about a situation where ``mathematics became the
mathematics of mathematized science of nature" (Cohen 1914
\cite[p.~22]{Co14}).  As Seidengart rightly observes,

\begin{quote}
\ldots independently of Kant, it was Leibniz who led Cohen along the
pathway of his ``logic of origin" (``Logik des Ursprungs"), which
\ldots is the logic of pure thought (``Logik des reinen
Denkens"). (Seidengart 2012 \cite[p.~131]{Se})
\end{quote}

However, Cohen's assessment of Leibniz had not always been thus
positive. Only for Cohen's later thought, from \emph{Logik} (1902)
onwards, did Leibniz's philosophy play a pre-eminent role.  In
contrast, for Cassirer, Leibniz had always been the philosophical hero
from the start of his philosophical career, as is evidenced by his
\emph{Leibniz' System}.  Eventually, however, Seidengart concludes,

\begin{quote}
\ldots inspite of the many innovations he was able to derive from
Leibniz's infinitesimal analysis, Cohen aligned in the end the
interpretation that young Cassirer laid out in his \emph{Leibniz'
System} and in his \emph{Erkenntnisproblem}, both of which were
explicitely cited by the founder of the Marburg School in 1914 [i.e.,
in (Cohen 1914 \cite[p.~24]{Co14})--the authors].
\end{quote}

Since Cohen (1914) was the last time that Cohen dealt with Leibniz's
infinitesimal analysis and its philosophical implications, this would
appear to suggest that Cohen and Cassirer agreed on matters Leibnizian
from the beginning of the 20th century until 1914 and perhaps even
later, until Cohen's death in 1918.

This, however, is not quite true as is already shown by Cohen's letter
of August 1910 where Cohen complained that Cassirer in SF had deviated
from the party line as he no longer recognized the primacy of the
concept of the infinitesimal.

Seidengart does not take into account the 1910 letter, and
concentrates on (Cohen 1914 \cite{Co14}).  However, a closer look at
Cohen (1914) reveals that this discrepancy of 1910, had \emph{not}
disappeared in 1914. True enough, in 1914 Cohen praised Cassirer's
\emph{Leibniz' System} (1902) and his \emph{Erkenntnisproblem} (1906)
as congenial elaborations of his own account of Leibniz and the role
of the infinitesimal.  More telling, however, is the fact that in 1914
Cohen did \emph{not} cite Cassirer's \emph{Substance and Function}
(Cassirer 1910)!  This omission suggests that the differences of 1910
between the two philosophers had not been resolved in the meantime.
Rather, Cohen implicitely recognized in 1914 that his interpretation
of Leibniz and that of Cassirer essentially differed.

Pursuing the opposite path, Skidelsky in his recent book \emph{Ernst
Cassirer.  The Last Philosopher of Culture} \cite{Ski} seeks to
emphasize the differences between Cassirer and Cohen, as well as the
alleged obsoleteness of the latter's infinitesimal account.  He too
invokes the 1910 letter, but leaves out the sentence in which Cohen
characterized SF as a possible substitute for the second part of his
\emph{Logik} (cf.~\cite[p.~64]{Ski}).  Skidelsky seeks to drive home
his case against Cohen by contending that ``Cohen's theory of
infinitesimals is in fact mistaken even from a purely mathematical
point of view, being based on an outmoded interpretation of calculus''
(\cite[p.~65]{Ski}).

The crown witness to Skidelsky's sweeping claim is, predictably,
Russell's \emph{The Principles of Mathematics}.  Describing an
infinitesimalist approach to the calculus as ``outmoded'' amounts
merely to toeing the line on the CDW approach to the formalisation of
analysis.  The alleged uniqueness of such an approach is being
increasingly challenged in the current literature.  See, e.g.,
(B\l{}aszczyk et al.~1013 \cite{BKS}); (Bair et al.~2013 \cite{B11}).

Cohen concluded his criticisms by urging Cassirer to ``take these
thoughts into intimate consideration in the new edition'' (of SF), but
such an edition never appeared.  Hence, Cassirer had to find another
opportunity where he could pay due respects to Cohen's philosophical
interpretation of infinitesimals.  He appears to have attempted to
overcome the clash with Cohen on matters infinitesimal in the article
\emph{Hermann Cohen und die Erneuerung der Kantischen Philosophie}
(Cassirer 1912 \cite{Cas12}), dedicated to Cohen on the occasion of
his 70th birthday.  In this paper he hailed Cohen as the innovator and
true heir of Kant's philosophy who had brought to the fore the
fundamental principle of mathematized natural science in terms of the
infinitesimal:

\begin{quotation}
Matter and movement, force and mass may be conceptualized in this
respect as instruments of knowledge. The high point of this
development is not, however, reached before we come back to the basic
mathematical motif underlying all specific conceptual formations of
the natural sciences. This motif presents itself to us in the
conceptual methodology of the ``infinitesimal". (Cassirer 1912
\cite[p.~260]{Cas12})
\end{quotation}

This appears to be a stronger endorsement of the infinitesimal
approach than it really was.  One should note that Cassirer spoke of
the ``methodology of the infinitesimal'' rather than the
``infinitesimal'' itself.  The methodology of the infinitesimal is
something more general than the infinitesimal itself, and Cassirer
seemed to have been well aware of this. Indeed, he sought to employ
this greater generality to bind Cohen's infinitesimal approach with
his own ``relational approach".  Actually he did not go beyond what he
had already offered in SF some years earlier, when he praised the
infinitesimal calculus as the first and most important example of the
many calculi developed in modern mathematics, e.g. Grassmann's
\em{Ausdehnungslehre}\em, Hamilton's theory of quaternions, the
projective calculus of distances, and many others (cf. SF, 95). From
this general perspective, then, it is easy to see that all these
calculi are ``relational'' or ``functional'' in a broad sense.
Meanwhile, the only example of a calculus that Cohen ever mentioned
was the infinitesimal calculus.  Superficially, Cassirer's remarks may
appear to be a reconciliation of the infinitesimal approach and the
relational approach, but they fail to convince.  For instance, Cohen
explicitly asserted that he was not primarily interested in the
infinitesimal calculus, but rather in the specific philosophical
ramifications of the concept of the infinitesimal which for Cohen
represented the triumph of pure thought (cf. \cite{Co83}, p. 32).
Hence, a general comparison of the infinitesimal \emph{calculus} with
other \emph{calculi} probably did not overly impress Cohen.  Cassirer
might have felt this inadequacy and offered a further argument:
 
 \begin{quotation}
Without [the mathematical leitmotif of the infinitesimal], it would
not be possible even to characterize rigorously the concept of
movement,%
\footnote{Strictly speaking Cassirer's claim is inaccurate.  It is
possible to get by in mathematics without infinitesimals, as CDW had
impressively shown.}
as it is presupposed by the mathematized natural science, to say
nothing of the task of fully comprehending the lawfulness of
movements.  Thereby the circle of the critical investigations is
closing.  Since without doubt, the concept of the infinitely small
does not denote a ``Being" that can be captured by the senses, but a
peculiar way and basic direction of thought:%
\footnote{This cryptic remark is meant to emphasize the non-empirical
character of the infinitesimal.  Unlike ideal realists as well as
empiricists, the Marburg school held that the infinitesimal was not to
be found ``out there'' in some empirical or ideal domain independent
of the cognizing subject; rather the infinitesimal was a way of
thinking or conceptualizing the world.}
but this basic direction is now revealed as the necessary
presupposition of the scientific object itself. (Cassirer 1912
\cite[p.~260-261]{Cas12})%
\footnote{Similar remarks appear already in (SF, 1910 (130ff), 1953
(99ff)).}
\end{quotation}

In fact, Cassirer never elaborated on the connection between Cohen's
``infinitesimal analysis" and his own ``relational analysis''.
Neither in SF nor in PSF did the concept of the infinitesimal play as
prominent a role as in Cohen's \emph{Prinzip} or his \emph{Logik der
reinen Erkenntnis}.  Rather, Cassirer used Cohen's ``methodology of
the infinitesimal'' only as a launch pad to develop his own
``methodology of the relational''.

Even more revealing is the fact that in the posthumously published
fourth volume of \emph{The Problem of Knowledge} (\cite{Cas50}),
written in the late thirties during his Swedish exile, he had
completely abandoned the infinitesimal standpoint; Cohen is not even
mentioned once.
%
%\footnote{It might be mentioned that by this time the pressure to
%conform to the mathematicians' contemporary views of foundations as
%having been by and large finalized by CDW had become enormous,
%particularly due to Russell's influence, and any philosopher of a
%lesser standing than Cassirer who challenged such a consensus risked
%ostracism.}
%

Cassirer's attempts in SF and elsewhere to connect his relational
account with the infinitesimal account of Cohen are to be judged as
less than fully convincing. Partisans of the infinitesimal approach
should not blame Cassirer for this shortcoming too harshly, however.
In his day, the effectiveness of an infinitesimal approach compared
with one based on epsilontics was not too compelling.  Moreover, the
conceptualization of infinitesimals as an idealizing completion was
not sufficiently understood to undertake a reasonable comparison with
other idealizing completions.  The advent of the various versions of
modern infinitesimal-enriched continua has changed the conceptual
landscape dramatically.

\section*{Acknowledgments}

T.~Mormann was partially funded by the research project FF2012-33550
of the Spanish government.  M.~Katz was partially funded by the Israel
Science Foundation grant no.~1517/12.  We are grateful to the
anonymous referees for a number of helpful suggestions.

\end{document}